\documentclass[a4paper,12pt,twoside]{article}
\usepackage[latin1]{inputenc}
\usepackage[T1]{fontenc}
\usepackage{amssymb}
\usepackage{dochier}
\usepackage[final]{lpretex}
\usepackage{title}
\usepackage{a4}
\usepackage{amscd}
\usepackage[all]{xy}

\EndOfDump

\def\DHpartformat#1{(#1) }
\def\DHrefpart#1{(\DHRefpart{#1})}
\LBsetstyleof{partno}{arabic}

\DocVersion[Moved out from Enriques ms; Sep  4 2006]{1}

\let\define\def

  \def\C {{\mathbb C}}

\def\GG {{\mathbb G}}   
  
  \def\P {{\mathbb P}}

\def\Z {{\mathbb Z}} 

\define \n {\mathbb N}
\define \z {\mathbb Z}
\define \q {\mathbb Q}
\define \PP {\mathbb P}

\def\sA {{\Cal A}}  \def\sC {{\Cal C}}
 \def\sE {{\Cal E}} \def\sF {{\Cal F}}
\def\sG {{\Cal G}} \def\sH {{\Cal H}} \def\sI {{\Cal I}}
 \def\sK {{\Cal K}} \def\sL {{\Cal L}}
\def\sM {{\Cal M}} \def\sN {{\Cal N}} \def\sO {{\Cal O}}
 
  \def\sU {{\Cal U}}
\def\sV {{\Cal V}} \def\sW {{\Cal W}} \def\sX {{\Cal X}}
\def\sY {{\Cal Y}}

\define \cN {\Cal N}
\define \cf {\Cal F}
\define \cg {\Cal G}
\define \cE {\Cal E}
\define \ce {\Cal E}
\define \cc {\Cal C}
\define \cV {\Cal V}
\define \cA {\Cal A}
\define \cK {\Cal K}
\define \cO {\Cal O}
\define \cF {\Cal F}
\define \cn {\Cal N}
\define \cI {\Cal I}
\define \sP {\Cal P}

\def\s {\sigma} \def\t {\theta}

\define \x {\xi}
\define \y {\eta}
\define \G {\Gamma}
\define \r {\rho}
\define \w {\omega}

\def\tX {\widetilde X}

\def \tC {\widetilde C}

\def \trho {\tilde {\rho}}

\def \tpi {\tilde{\pi}}

\def \tp {\widetilde{\mathbb P}}
\define \tH {\widetilde H}
\define \tG {\widetilde{\Gamma}}
\define \tW {\widetilde W}
\define \tF {\widetilde F}
\define \tm {\tilde m}
\define \St {\widetilde S}
\define \Xt {\widetilde X}
\define \tS {\widetilde S}
\define \tpsi {\tilde \psi}
\define \tL {\widetilde L}
\define \tE {\widetilde E}
\define \tl {\tilde l}
\define \tA {\widetilde A}
\define \tom {\tilde\omega}
\define \tT {\widetilde T}
\define \tB {\widetilde B}
\define \tf {\widetilde f}
\define \tsA {\widetilde{\sA}}
\define \tsigma{\widetilde{\sigma}}
\define \tM {\widetilde M}
\define \tphi {\widetilde{\phi}}
\define \trho {\widetilde{\rho}}
\define \tR {\widetilde R}
\define \tp {\tilde p}
\define \tq {\tilde q}
\define \tc {\tilde c}
\define \tsF {\widetilde {\sF}}
\define \tx {\tilde x}
\define \tg {\tilde g}
\define \tw {\tilde w}
\define \ts {\tilde s}
\define \tdelta {\widetilde\delta}
\define \talpha {\widetilde\alpha}
\define \tbeta {\widetilde\beta}
\define \tsG {\widetilde{\sG}}
\define \tth {\widetilde{\Theta}}
\define \tchi {\widetilde{\chi}}

\def\pd {\partial}

\def \Dx1 {\frac{\pd}{{\pd} x_1}}
\def \Dy1 {\frac{\pd}{{\pd} y_1}}
\def \Dz1 {\frac{\pd}{{\pd} z_1}}
\def \Dx2 {\frac{\pd}{{\pd} x_2}}
\def \Dy2 {\frac{\pd}{{\pd} y_2}}
\def \Dz2 {\frac{\pd}{{\pd} z_2}}

\def\q {\quad}

\def\mapdiagr#1{\Big\searrow\rlap{$\raise 5pt\vbox{{\hbox{$\mkern -15mu\scriptstyle#1$}}}$}}   

\def\mapdiagl#1{\llap{$\raise 5pt\vbox{{\hbox{$\scriptstyle#1\mkern
-15mu$}}}$}\Big\swarrow}              

\def\Mapdiagr#1{\nearrow\rlap{$\lower 5pt\vbox{{\hbox{$\mkern
-15mu\scriptstyle#1$}}}$}} 

\def\Mapdiagl#1{\llap{$\lower 5pt\vbox{{\hbox{$\scriptstyle#1\mkern
-15mu$}}}$}\searrow} 

\def\Mapswr#1{\swarrow\rlap{$\lower 5pt\vbox{{\hbox{$\mkern
-15mu\scriptstyle#1$}}}$}}              

\def\Mapnwl#1{\nwarrow\rlap{$\lower 5pt\vbox{{\hbox{$\mkern
-15mu\scriptstyle#1$}}}$}}

\def \inj {\hookrightarrow}

\define \Rhook {\hookrightarrow}

\def \half {\raise1pt\hbox{$\scriptstyle
        \frac{1}{2}\displaystyle$}}

\def \x{{\sl X}\llap{$\mkern -2mu {\scriptstyle -}$}}
\def \End {\operatorname{End}}
\def \Hom {\operatorname{Hom}}

\def \Symm {\operatorname{Sym}}
\def \Res {\operatorname{Res}}

\def \Pic {\operatorname{Pic}}

\define \Kod {\operatorname{Kod}}
\define \dimension {\operatorname{dim}}
\define \codim {\operatorname{codim}}
\define \contr {\operatorname{contr}}
\define \rk {\operatorname{rank}}
\define \im {\operatorname{im}}
\define \Mor {\operatorname{Mor}}
\define \Cl {\operatorname{Cl}}
\define \Hilb {\operatorname{Hilb}}
\define \degree {\operatorname{deg}}
\define \mult {\operatorname{mult}}
\define \Aut {\operatorname{Aut}}
\define \NS {\operatorname{NS}}
\define \Gal {\operatorname{Gal}}
\define \ch {\operatorname{char}}
\define \Jac {\operatorname{Jac}}
\define \Km {\operatorname{Km}}
\define \Sec {\operatorname{Sec}}
\define \Stab {\operatorname{Stab}}
\define \Br {\operatorname{Br}}
\define \inv {\operatorname{inv}}
\define \tr {\operatorname{tr}}
\define \Frob {\operatorname{Frob}}
\define \Symn {\operatorname{Sym}^n}
\define \Ev {\sE^\vee}
\define \ordp {\operatorname{ord}_p}
\define \Supp {\operatorname{Supp}}
\define \Ann {\operatorname{Ann}}
\define \disc {\operatorname{disc}}
\define \Lie {\operatorname{Lie}}
\define \embdim {\operatorname{embdim}}

\def\Ext{\operatorname{Ext}}

\def\Chow{\operatorname{Chow}}
\def\Nm{\operatorname{Nm}}
\def\barC{\overline{C}}

\def\hod#1#2#3#4{\ensuremath{\if#30 H^{#2}({#1},{\cal O}_{#1}) \else 
 H^{#2}(#1,\Omega^{#3}\if\relax{#4}\relax_{#1}\else _{#1/#4}\fi)\fi}}
\begin{document}
\title[Thomae's formulae]{Thomae's formulae for non-hyperelliptic curves and 
spinorial square roots of theta-constants on the moduli space of curves}
\author{N. I. Shepherd-Barron}
\address{D.P.M.M.S.\\
Cambridge University\\
Cambridge CB3 0WB\\
U.K.}
\email{nisb@dpmms.cam.ac.uk}

\maketitle
\begin{section}{Introduction}\label{intro}
\medskip
One aim of this paper
is to give formulae for the ratios of
the $8$th powers of the $2$nd order theta functions,
and so the theta-nulls,
associated to the Jacobian of a genus $g$ curve $C$ over $\C$
in terms of determinants of twisted pluricanonical differentials;
see \ref{thetanull} below. In fact,
we give formulae for the $2N^2$ powers of the $N$th order
theta functions for any even integer $N$.
More generally, we give formulae in terms of such determinants
for ratios of algebraic versions of theta functions 
over a field $k$ of characteristic $p$ that is prime to $N$,
for $N$ even or odd. 

In order to be able to substitute effectively into these formulae,
starting, say, from the equations of the canonical model of 
a non-hyperelliptic curve $C$, we must know also the $N$-torsion
points $P$ on $\Jac^0C$. It turns out that we can,
starting from the equations of $C$, write down the equations
that define multiplication by $N$ on $\Jac^0_C$, and so determine the 
$N$-torsion points. We can also derive
equations for the group law on $\Jac^0C$ in projective terms,
starting from the equations of $C$, but this is more complicated.

These formulae provide a polynomial
formula for the Torelli morphism $\sM_g\to A_g$, where $\sM_g$
is the stack of smooth curves of genus $g$ and $A_g$ is the coarse moduli
space of principally polarized abelian varieties (more accurately,
they do so at level $(N,2N)$ if $N$ is even and at level $N$
if $N$ is odd). This is because the $N$th
order theta constants give canonical homogeneous co-ordinates
on $A_{g,(N,2N)}$ or $A_{g,N}$; we recall below how to do this
in a way that includes characteristic $2$. 
For hyperelliptic curves, and when $N=2$,
this is what Thomae's formulae do, although his formulae 
also give, in transcendental terms, 
the common scalar factor that our approach neglects.

The idea of writing theta functions on $\Jac^{g-1}C$
at points of $C^{(2n-1)(g-1)}$ in terms of determinants of $n$-fold
differentials for $n\ge 2$ is taken from
Proposition 4.1 of Matone and Volpato's paper \cite{MaV}. 
That is a specialization
of an addition theorem due to Fay (Proposition 2.16 of \cite{F}).
This, in turn, generalizes a formula due to Klein 
(\cite{K}, vol. III, p. 429),
who takes $N=2$ and describes theta functions on $\Jac^{\rho(2g-2)}C$, or
$\Jac^0C$,
in terms of $\rho(2g-2)\times \rho(2g-2)$ determinants
in $H^0(C,\sO(\rho K_C+L))$, where $L$ is a theta-characteristic
and $\rho$ a positive integer. 

Klein goes on to describe a relation between the theta-constants
of a plane quartic and its discriminant. 
This part of his paper was largely an 
investigation of the behaviour of theta functions
as a genus $3$ curve degenerates, but it can also be read as the
start of an attempt to find formulae for the theta-constants of a genus $g$
curve in terms of its projective invariants rather than in terms of the equations
of its canonical model with respect to some arbitrary choice of basis
in $H^0(C,\omega_C)$. For example, one might consider 
configurations of higher Weierstrass points, which are
known to determine the moduli of the curve (cf. \cite{GIT})
and then seek formulae
for the theta-constants in terms of their invariants;
in fact, such formulae
would be truer analogues of Thomae's, since his involve the cross-ratios
of the branch points determined by the hyperelliptic involution.
On the other hand, this paper obviates the task of computing invariants
of Weierstrass points from the equations of the curve
and second, since there are only $N^{2g}$
theta-constants to compute, it might be that they (with the
Torelli theorem) are better adapted than is invariant theory
to the task of detecting whether two curves given by equations are isomorphic.

What we do here is to twist differentials in a slightly different way
and to make things more algebraic. We do this by considering
ratios of theta functions rather than the thetas themselves.
This means that our results hold over any field. By comparison, the formulae 
of Klein \emph{et al.} involve the prime form $E$,
some exponential factors that Fay denotes by $\sigma$
and some constants that depend, in an undetermined way,
on the choice of basis of vector spaces $H^0(C,\sL)$
for various invertible sheaves $\sL$ on $C$.

The prime form has a purely algebraic construction
in terms of a theta characteristic of multiplicity exactly $1$,
and so makes sense in all characteristics except $2$,
where such theta characteristics do not always exist.
The exponential factors are explicit, and have neither zeroes nor poles,
but I do not know what they mean algebraically.
Taking ratios has the advantage of
removing that difficulty, but it also loses the prime forms.
(Precisely, $\sigma$ is defined in terms of the prime form $E$ 
and normalized differentials $\omega_i$ by
$\sigma(z)=\exp(-\sum \int_{A_i}\omega_i(z)E(z,w)).$ Of course, the
process of normalizing the differentials is also transcendental.)
To pin down the constants we use Mumford's theta groups.

There are \cite{Sz} other well known determinantal descriptions of the theta
divisor $\Theta$ on a Jacobian $J$, which amount to giving a homomorphism
of locally free sheaves on $J$ whose degeneracy
locus is $\Theta$. As is also well known,
these descriptions extend to stacks of $G$-bundles
on a curve $C$ for other reductive groups $G$. 
In particular, when Mumford's approach to theta-characteristics
\cite{MuTh} is globalized in this way, it leads to the second main
point of this paper: Tsuyumine showed (\cite{Ts}, Theorem $1$) that theta-constants
on the moduli space $M_{g,(4,8)}$ of curves with level $(4,8)$ structure,
which are determinants,
have a square root, and we show here that these
square roots have a local description as pure spinors. 

The outline of our approach is this. Assume that
$g\ge 2$; our formulae make no sense when $g=1$. 
Given an $N$-torsion point $P$ on $A:=\Jac^0C$, 
we write down,
in terms of determinants of twisted quadratic differentials,
an explicit rational function $f_P$
on $X:=\Jac^{g-1}C$
with an equality $(f_P)=N\Theta_P-N\Theta_0$ of divisors, where
$\Theta_0=\Theta$ is the theta divisor on $X$ and $\Theta_P=t^*_P\Theta_0$
is its translate by $P$. These functions, defined so far up to
an arbitrary scalar, give a formula for the Weil pairing
$e_N$. Then, assuming that we know how to write down all the $N$-torsion points
of $A$, we choose {\emph{normalized}} functions ${\tf}_P$, in a way that 
involves Mumford's theta groups.  The $N^{2g}$ functions ${\tf}_P$, with ${\tf}_0=1$, 
define a morphism (an embedding if $N\ge 3$)
$X\to \P^{N^{2g}-1}$ to a projective space with a co-ordinate system
that has been specified by the theta group structure.
The linear span of the image of $X$ is a copy of $\P^{N^g-1}$,
which varies with the moduli of $C$, or $A$.

To make this work we need to be able to find the $N$-torsion points $A[N]$.
We do this in two ways. First, we give a projective
description of the group law on $\Jac_C^0$, which in turn
is given in terms of a chord and tangent construction on $\Jac_C^{g-1}$,
and from this derive equations defining $A[N]$. Second, which is more efficient,
we use interpolation to find the functions that define multiplication by $N$,
and then derive equations of $A[N]$.

Now assume that $p\ne 2$. 
Then the moduli point is obtained by evaluating the functions ${\tf}_P$
at the theta characteristics, which are the fixed points on $X$
of the involution provided by $[-1_A]$, acting as $D\mapsto K_C-D$.
Over $\C$, these values are just the ratios of the $N$th powers
of the theta-constants $\theta[a,b](p+\tau q,\tau)$,
where $a,b$ run over $((1/N)\Z/\Z)^g$ and $p,q$ over
$(1/2)\Z/\Z)^g$. When $p=2$ we do something slightly different,
described below.

I am very grateful to Robin de Jong, for telling me at the Schiermonnikoog
conference about Matone and Volpato's 
paper \cite{MaV}. I am also very grateful to GNSAGA at the University of Rome, La Sapienza,
and the Institut Mittag-Leffler (Djursholm) for their support.

\end{section}
\bigskip
\begin{section}{Principal symmetric abelian torsors}\label{prelim}
\medskip
This section is mostly a summary of well known facts. The only thing in
it that might be new is to extend the use of theta-constants to co-ordinatize
suitable level covers of $\sA_g$ to include characteristic $2$, that is,
to cover all of $\Sp \Z$ rather than $\Sp \Z[1/2]$.

We use Alexeev's idea of identifying the stack of principally polarized
abelian schemes with the stack of torsors with a suitable specified
divisor (\cite{Al}, Corollary $3.0.7$). We impose the additional constraint
of symmetry.

An {\emph{abelian torsor}} is a torsor under an abelian scheme $A\to S$.
Given an abelian scheme $A\to S$, an $A$-torsor $X\to S$ (assumed projective)
is {\emph{symmetric }} if the action
of $A$ on $X$ is provided with an extension to an action of the semi-direct product
$A\rtimes [-1_A]$, where $[-1_A]$ is the involution $x\mapsto -x$
of $A$. Given such a torsor, we let $Fix_X$ denote the fixed-point subscheme
of $X$ with respect to $[-1_A]$; it is a torsor under the $2$-torsion
subscheme $A[2]$ of $A$. 

Note that a smooth projective torsor $X\to S$ determines the abelian scheme
$A$ under which it is a torsor, by $A=\Aut^0_{X/S}$, the connected component of
the automorphism group scheme.
\def\underscore{\underline}
\begin{lemma} if $X\to S$ is an abelian torsor, then the sheaf 
${\underscore{\Pic}}^\tau_X$ is 
represented by the dual abelian scheme $\Pic^0_A$.
\begin{proof} The action of $A$ on $X$ is an isomorphism $A\times_S X\to X\times_S X$
that commutes with the second projections to $X$.
This defines an isomorphism of abelian sheaves 
$${\underscore{\Pic}}^\tau_X\times
{\underscore{\Pic}}^\tau_X \to {\underscore{\Pic}}^\tau_A\times {\underscore{\Pic}}^\tau_X;$$
taking cokernels by $pr_2^*{\underscore{\Pic}}^\tau_X$ gives the result,
since ${\underscore{\Pic}}^\tau_A$ is represented by $\Pic^0_A$.
\end{proof}
\end{lemma}

\begin{lemma} \label{polarize}
An ample line bundle $\sM$ on an abelian torsor $X\to S$ defines
a polarization $\lambda$ on $A:=\Aut^0_{X/S}$.
\begin{proof} Define $\lambda:A\to {\underscore{\Pic}}^\tau_X$
by $\lambda(a)=t_a^*\sM\otimes\sM^{-1}$, where $t_a$ denotes ``translation by $a$''.
The previous lemma completes the proof.
\end{proof}
\end{lemma}

\begin{definition} An effective divisor $\Theta$ on the abelian torsor
$X\to S$ is \emph{principal} if it defines a principal polarization
on $A:=\Aut^0_{X/S}$. A \emph{principal symmetric abelian torsor}
is a symmetric abelian torsor $X\to S$ with an effective principal divisor
$\Theta$ that is symmetric, that is, preserved by the action of $[-1_A]$.
\end{definition}

Let $\sX_g$ be the stack of principal symmetric abelian torsors of dimension
$g$, $\sA_g$ the stack of principally polarized abelian varieties 
(ppav's) of dimension $g$ and $\sU\to\sA$ the universal ppav.
By \ref{polarize}, there will be a natural morphism $\beta:\sX\to\sA$
once it has been defined on morphisms. Of course, if $f:(X,\Theta)\to (Y,\Phi)$
is a morphism in $\sX_g$, then $g=\beta(f):\Aut^0_X\to\Aut^0_Y$ 
is defined at the level of points by $g(a)(y)=f(a(f^{-1}(y)))$.

Recall (\cite{EHSB}, Lemma 1.6) that, for any morphism 
$F:\sY\to\sX$ of algebraic stacks,
there is a (representable) groupstack $\sK_F\to\sY$, the inertia stack or
automorphism stack or kernel, such that
for any space $x:U\to\sY$, the product $\sK_F\times_{\sY,x}U\to U$
is the automorphism
group scheme of the object $x$ over $F(x)$. If $F$ is a gerbe
(\cite{EHSB}, Definition 1.10) and $\sK_F$ is trivial,
then $F$ is an isomorphism.

The next result is analogous to Corollary $3.0.7$ of \cite{Al};
the only difference is that we take the torsors to be symmetric.

\begin{proposition}\label{isom} $\beta:\sX_g\to\sA_g$ is an isomorphism. 

\begin{proof} We show first that 
$\beta$ is a gerbe.

To prove local surjectivity on objects, suppose that
$(A,\lambda)\to S$ is a principally polarized abelian scheme.
There is a faithfully flat map $T\to S$ such that on
$A_T$ there is a symmetric line bundle $\sO(\Theta)$
which defines $\lambda_T$. Then $(A_T,\Theta)$ provides
a suitable object in $\sX$.

To prove local surjectivity on morphisms, suppose that
$(X,\Theta)\to S$ and $(Y,\Phi)\to S$ are objects in $\sX$
and $\alpha:(A=\Aut^0_X,\lambda)\to(B=\Aut^0_Y,\mu)$ is an isomorphism of
the corresponding principally polarized abelian schemes.
Then, locally on $S$, 
there are isomorphisms $\xi:A\to X$ and $\eta:B\to Y$
such that $\lambda$ is defined by $\xi^*\Theta$ and
$\mu$ by $\eta^*\Phi$. Then $\zeta:=\eta\circ\alpha\circ\xi^{-1}:X\to Y$
is an isomorphism such that $\zeta^{-1}(\Phi)$ and $\Theta$
define the same polarization. Then they differ by a translation;
that is, there is an $S$-point $a$ of $A$ 
with $\zeta^{-1}(\Phi)=a^*\Theta$, so that $\zeta\circ a^{-1}$
is an isomorphism $\psi:X\to Y$ with $\psi^{-1}(\Phi)=\Theta$.

So $\beta$ is a gerbe, and it is enough to check that
$\beta$ is locally injective on automorphisms. For this,
we can assume that $S$ is a geometric point.

Suppose that $\xi:(X,\Theta)\to(X,\Theta)$ is an automorphism
that induces the identity on $A$. That is,
$a(x)=\xi a\xi(x)$ for all $x\in X$ and $a\in A$.
So $\xi a=a\xi$ as automorphisms of $(X,\Theta)$
for all $a\in A$. Locally on $S$, there is an isomorphism
$X\to A$; via this, $\xi$ is an isomorphism
$\xi:A\to A$ (of varieties, not of
abelian varieties) such that $\xi\circ t_a=t_a\circ\xi$ for all
$a$, where $t_a$ denotes translation by $a$.
So $\xi(a+b)=a+\xi(b)$ for all $a,b\in A$.
Now there is an automorphism
$\gamma$ of $A$, as an abelian scheme, and a point $c_0$ on $A$
such that $\xi(x)=\gamma(x)+c_0$ for all $x$.
Since $\gamma(a+b)=\gamma(a)+\gamma(b)$, substituting into
$\xi(a+b)=a+\xi(b)$ gives $\gamma(a)=a$ for all $a$.
\end{proof}
\end{proposition}

There is usually no line bundle on a ppav $(A,\lambda)$ that defines
$\lambda$. However, there is a unique totally symmetric line bundle
$\sM_2$ that defines $2\lambda$, namely $\sM_2=\sL^\Delta(\lambda)$.
Note that $2\lambda$ determines $\lambda$,
because $\End(A)$ is torsion free.
This leads to consideration of the $2$-commutative
diagram
$$\xymatrix{
{\sX_g}\ar[r]^{\beta}\ar[d]_{\pi} & {\sA_g}\ar[d]^{\alpha}\\
{{}_2\sX_g} \ar[r]^{\delta} &{{}_2\sA_g.}
}$$
Here ${}_2\sX_g$ is the stack associated to the prestack of pairs $(X,\sL_2)$,
where $X\to S$ is a symmetric abelian torsor and $\sL_2$ is a line bundle on
$X$, trivial along $Fix_X$, taken up to isomorphism
and taken modulo the pullback
of line bundles on $S$, that defines twice a principal polarization
$\lambda$ on $A:=\Aut^0_X$, and ${}_2\sA_g$ is the stack of
ppav $(A,\lambda)$ with a totally symmetric line bundle $\sM_2$
as above. It is immediate that $\alpha$ is an isomorphism,
so we identify ${}_2\sA_g$ with $\sA_g$.

\begin{proposition} ${}_2\sX_g$ is algebraic and
$\delta:{}_2\sX_g\to\sA_g$ makes
${}_2\sX_g$ isomorphic to the classifying stack $B\sU[2]$,
where $\sU\to\sA_g$ is the universal abelian scheme.
\begin{proof} First, the obvious forgetful map
provides a section of $\delta$.
So it is surjective on objects.

The proof that $\delta$ is locally
surjective on morphisms is exactly the same as in \ref{isom},
except for replacing $\Theta$ by $\sL_2$ and $\Phi$ by $\sM_2$. 

So $\delta$ is a neutral gerbe. It is clear that the automorphism
group of $(X,\sL_2)$ that covers the identity on $A=\Aut^0_X$ is $H(\sL_2)$,
which is $A[2]$, so that (\cite{EHSB}, Prop. 1.11 or \cite{LMB}, Lemme 3.21)
${}_2\sX\cong B\sU[2]$. It follows that ${}_2\sX_g$ is algebraic. 
\end{proof}
\end{proposition}
\end{section}
\bigskip
\begin{section}{Level structures and moduli}\label{level structures}
\medskip
\begin{subsection}{Abstract theta groups}\label{abstract theta groups}
\smallskip
We now review Mumford's theory of theta groups,
with slight modifications imposed by the point of view that line bundles,
divisors, rational functions and linear systems live on 
$X$ rather than on $A$.
$N$ will always denote an integer with $N\ge 2$ that is invertible
in whatever base ring, or scheme, we are considering.

Put $L=(\Z/N\Z)^g\times\mu_N^g$ with its standard symplectic pairing
$e:L\times L\to\mu_N$. As usual, a \emph{level $N$ structure} on $A$
is a choice of isomorphism $\phi:L\to A[N]$ taking $e$ to
the Weil pairing $e_N$;
a level $N$ structure on $X$ is a level $N$ structure on $A=\Aut^0_X$.

There are bilinear pairings $d:L\times L\to\mu_N$
whose skew-symmetrization is $e$, that is, that satisfy 
$e(P,Q)=d(P,Q)/d(Q,P)$; for example, choose a primitive $N$th
root $\zeta$ of unity and a symplectic
basis $\{P_1,...,P_g,Q_1,...,Q_g\}$ of $L$ and then define $d$
by $d(P_i,Q_i)=\zeta, d(Q_i,P_i)=1$ and $d=1$ on all other
pairs of basis elements. Of course, if $N$ is odd,
then there is a canonical choice of such a form $d$,
namely, $d(P,Q)=e(P,Q)^{(N+1)/2}$. 

Mumford \cite{MuEqnsI} defines the (abstract) theta group $\sG_L$,
a central extension of $L$
by $\GG_m$, to be $\sG_L=\GG_m\times L$, with group law
$$(\lambda,P).(\nu,Q)=(\lambda\nu\ d(P,Q),P+Q)$$
and natural surjection $\pi:\sG_L\to L$.
The commutator pairing is $e$. The next lemma shows that,
given $e$, $\sG_L$ is independent of the choice of $d$.

\begin{lemma}\label{theta group uniqueness}
Given a symplectic pairing $e$ on $L$, there is 
a unique central extension of $L$ by $\GG_m$ whose commutator
pairing is $e$.
\begin{proof} 
By \cite{Br}, pp. $97$ and $127$, there is a short exact sequence
$$0\to \Ext^1_{comm.\ gp.\ sch.}(A[N],\GG_m)\to H^2(A[N],\GG_m)
\to\Hom(\bigwedge{{}^2}A[N],\GG_m)\to 0.$$ The middle group classifies central 
extensions of $A[N]$ by $\GG_m$ and the left hand group is trivial 
(\cite{MuAV}, p. $223$, Lemma $1$).
\end{proof}
\end{lemma}

Such a central extension is a {\emph{theta group}}.

Say that an automorphism of $\sG_L$ is {\emph{quasi-trivial}}
if it induces the identity on the central $\GG_m$ and on the quotient $L$.
As usual, write $\Hom(L,\GG_m)=L^\vee$.

\begin{lemma}\label{theta group automorphism}
\part[i] If $\chi\in L^\vee$, then there is a quasi-trivial
automorphism $\alpha=\alpha_\chi$ of $\sG_L$ defined by
$\alpha(\lambda,\ P)=(\lambda\chi(P),\ P)$.

\part[ii] The map $\chi\mapsto\alpha_\chi$ is an isomorphism from
$L^\vee$ to the group of quasi-trivial automorphisms of $\sG_L$.
\begin{proof} \DHrefpart{i} is trivial.

For \DHrefpart{ii}, suppose that $\alpha$ is a quasi-trivial automorphism of $\sG_L$.
Then $\alpha(\lambda,\ P)=(\lambda\beta(\lambda,P),\ P)$ for some morphism
$\beta:\GG_m\times L\to\GG_m$ of schemes. Write $\beta(\lambda,P)=\beta_P(\lambda)$ and
$\GG_m=\Sp k[t,t^{-1}]$; then
$\beta_P^*(t)=\rho t^n$ for some $\rho\in k^*$ and $n\in\Z$.
That is, $\beta(\lambda,P)=\rho(P)\lambda^{n(P)}$
for some maps $\rho:L\to\GG_m$ and $n:L\to\Z$.

Since $\alpha$ is a homomorphism,
$$\beta(\lambda\mu d(P,Q),P+Q)=\beta(\lambda,P)\beta(\mu,Q).$$
This gives
$$\rho(P+Q)(\lambda\mu d(P,Q))^{n(P+Q)}=\rho(P)\rho(Q)\lambda^{n(P)}\mu^{n(Q)}$$
for all $\lambda,\mu\in\GG_m$ and $P,Q\in L$. By fixing $\mu,P$ and $Q$
and letting $\lambda$ vary, we see that
$n$ is constant, so that
$$d(P,Q)^{-n}=\rho(P+Q)\rho(P)^{-1}\rho(Q)^{-1}.$$
In particular, $d^n$ is symmetric, so that
$$\zeta^n=d(P_i,Q_i)^n=d(Q_i,P_i)^n=1.$$
So $N$ divides $n$ and $d^n=1$. Thus $\rho:L\to\GG_m$
is a character and $\beta(\lambda,\ P)=\lambda^n\rho(P)$.
Since $\beta(\lambda,0)=1$, it follows that $n=0$ and we are done.
\end{proof}
\end{lemma}
\end{subsection}
\medskip
\begin{subsection}{Concrete theta groups}\label{concrete theta groups}
Suppose that $(X,\Theta)$ is a principal symmetric abelian torsor
with $A=\Aut^0_X$. For $P\in A$, define $t_P:X\to X$ to be 
the translation by $P$ (or action by $P$) and $\Theta_P$ to be the translate
$\Theta_P=t_P^*\Theta$. In particular, $\Theta_0=\Theta$. Then
the {\emph{level $N$ theta group of}} $(X,\Theta)$, to be denoted by $\sG_{X,\Theta,N}$,
is the set
of pairs $(f_P,P)$, where $P\in A[N]$ and $f_P$
is an isomorphism $t_P^*\sO_X(N\Theta)\to\sO_X(N\Theta)$.
Over a field, such an $f_P$ can and will be regarded as
a rational function on $X$ whose divisor $(f_P)$ is
$(f_P)=N\Theta_P-N\Theta_0$, where $\Theta_0=\Theta$.
Such a function is a {\emph{Weil function for}} $P$ (at level $N$).
The kernel of the projection $\sG_{X,\Theta,N}\to A[N]$
is naturally identified with $\GG_m$.
The group structure on $\sG_{X,\Theta,N}$ is defined by
$(f_P,P)(f_Q,Q)=(f_P.t_P^*f_Q,P+Q)$,
$(f_P,P)^{-1}= (t_{-P}^*f_P^{-1},-P)$ and by taking $(1,0)$
as the identity. It is clear that $\sG_{X,\Theta,N}$ is a theta group.

\begin{lemma}\label{ext}\label{weil pairing}
\part[i] The Weil pairing $e_N$ determines and is determined by
the group scheme $\sG_{X,\Theta,N}$, when this is regarded
as a central extension of $A[N]$ by $\GG_m$.

\part[ii] The Weil pairing $e_N$ is given by
$e_N(P,Q)=\frac{f_P.t_P^*f_Q}{f_Q.t_Q^*f_P}$.
\begin{proof} \DHrefpart{i} is just a restatement of 
\ref{theta group uniqueness}.

For \DHrefpart{ii}, according to \cite{MuAV}, p. $227$,
the commutator pairing on $\sG_{X,\Theta,N}$ determines $e_N$.
\end{proof}
\end{lemma}

\begin{definition} Given a principal symmetric abelian torsor $(X,\Theta)$, or a ppav $(A,\lambda)$,
a {\emph{level $(N,\vartheta)$ structure}} on $(X,\Theta)$ is a choice
of isomorphism $\psi:\sG_L\to\sG_{X,\Theta,N}$ that induces
the identity on the central $\GG_m$s.
We denote by $\sA_{g,N,\vartheta}$ the stack of $g$-dimensional 
ppav's (or principal symmetric abelian torsors)
with a level $(N,\vartheta)$ structure and, as usual, by $\sA_{g,N}$
the stack of $g$-dimensional principal symmetric abelian torsors (or ppav's) with a level $N$ structure.
\end{definition}

\begin{proposition} The character group $L^\vee$ acts freely on $\sA_{g,N,\vartheta}$ 
and the quotient is identified with $\sA_{g,N}$.
\begin{proof} This follows at once from
\ref{theta group automorphism}.
\end{proof}
\end{proposition}

\begin{definition}\label{Weil set} 
Suppose that $(X,\Theta)$ is a principal symmetric abelian torsor and that, as usual, $A=\Aut^0_X$.

\part[i] A {\emph{Weil set}}
is an ordered set $(f_P)_{P\in A[N]}$
of Weil functions. 

Now suppose fixed a level $(N,\vartheta)$ structure on $(X,\Theta)$.
Fix also a bilinear pairing $d$ on $A[N]$ with $e_N(P,Q)=d(P,Q)/d(Q,P)$.

\part[ii] The Weil set
$(\tf_P)_{P\in A[N]}$ is {\emph{normal}} if
$\tf_P.t_P^*\tf_Q=d(P,Q)\tf_{P+Q}$
and $\tf_{-P}=t_{-P}^*\tf_P^{-1}$
for all $P,Q\in A[N]$.

\part[iii] A {\emph{normalization}} of a given Weil set
$(f_P)_{P\in A[N]}$
is an ordered set of non-zero scalars $(\alpha_P)_{P\in A[N]}$
such that $(\alpha_Pf_P)$ is normal.
\end{definition}

Note that for any Weil set $(f_P)_{P\in A[N]}$,
each ratio $f_P.t_P^*f_Q/f_{P+Q}$ has trivial divisor, 
and so is a constant function.

Given a Weil set, there is, by \ref{theta group uniqueness}, 
an isomorphism of schemes
$\psi:\sG_L\to\sG_{X,\Theta,N}$
defined by $(\lambda, P)\mapsto (\lambda f_{\phi(P)},\phi(P)).$
Notice that $\psi$ is a 
level $(N,\vartheta)$ structure if and only if $(\tf_P)_{P\in A[N]}$
is a normal Weil set; this is the reason for defining the latter notion.

In particular, there exists an 
isomorphism
$\psi:\sG_L\to\sG_{X,\Theta,N}$, and so there
exists a normalization $(\alpha_P)$ of any given Weil set $\{f_P\}$.
To calculate one, we proceed as follows.

By associativity, the conditions on the $\alpha_P$ are that
$$\frac{\alpha_P.\alpha_Q}{\alpha_{P+Q}}=\frac{f_{P+Q}}{f_P.t_P^*f_Q}.d(P,Q)$$
for all $P,Q\in A[N]$.

Note that $\tf_0=1$ automatically.
Start by imposing $f_0=1$ and $\alpha_0=1$.
Choose a point $x_0$ on $X$ in general position.

\begin{lemma}\label{powers}
$$\tf_P^N = \epsilon(P)\frac{\prod_{m=1}^{N-1}f_{mP}(x_0)}{f_P^N(x_0)
\prod_{m=1}^{N-1} f_{mP}(x_0-P)}f_P^N,$$
where $\epsilon(P)=1$ if $N$ is odd and $\epsilon(P)=\pm 1$ if $N$ is even. 

\begin{proof} All ratios $\frac{\alpha_P.\alpha_Q}{\alpha_{P+Q}}$ 
are determined by evaluating 
$\frac{f_{P+Q}}{f_P.t_P^*f_Q}$ at $x_0$ and evaluating $d(P,Q)$; 
taking $Q$ to be the various multiples of $P$
and multiplying the products 
$\frac{\alpha_P^2}{\alpha_{2P}},\frac{\alpha_P.\alpha_{2P}}{\alpha_{3P}},
\ldots, \frac{\alpha_P.\alpha_{(N-1)P}}{1}$ then gives 
$$\tf_P^N = \frac{\prod_{m=1}^{N-1}f_{mP}(x_0)}{f_P^N(x_0)
\prod_{m=1}^{N-1} f_{mP}(x_0-P)}d(P,P)^{N(N-1)/2}f_P^N,$$
from which the result follows.
\end{proof}
\end{lemma}

This formula also shows that to compute the functions
$\tf_P$, we can assume that each $\alpha_P$ lies in $\mu_N$.
To find these $N$th roots is, by \ref{theta group automorphism}, a matter
of splitting a torsor under $L^\vee$. This can be done explicitly by
adapting a little of Igusa's treatment of
theta characteristics (\cite{I}, top of p. 210).

Define $\gamma(P,Q)=d(P,Q)\frac{f_{P+Q}}{f_Pt_P^*f_Q}$ for $P,Q\in A[N]$.
We must solve the equations
$$\frac{\alpha_P\alpha_Q}{\alpha_{P+Q}}=\gamma(P,Q)$$
for all $P,Q\in A[N]$.

Write $\alpha(P)$ in place of $\alpha_P$ and put $P_i=r_i$, $Q_i=r_{i+g}$
for $i=1,\ldots,g$.
For each $i=1,\ldots,2g$, assign $\alpha(r_i)\in\mu_N$
arbitrarily. Then for any integer $m$,
$$\alpha((m+1)r_i)=\alpha(mr_i).\alpha(r_i)/\gamma(mr_i,r_i).$$
So $\alpha(mr_i)$
is determined. Next, if $2\le q\le 2g$, then
$$\alpha(\sum_1^q m_{i_j}r_{i_j})=\alpha(\sum_1^{q-1}m_{i_j}r_{i_j}).\alpha(m_{i_q}r_{i_q})/
\gamma(\sum_1^{q-1}m_{i_j}r_{i_j},m_{i_q}r_{i_q}),$$ 
so that, by induction on $q$, $\alpha$ is determined as a function $\alpha:A[N]\to\mu_N$.

It is convenient to have also the notion of a {\emph{symmetric}}
level $(N,\vartheta)$ structure. 
For even $N$ this is described on pp. $194-6$
of \cite{GIT} in the context of abelian varieties, and we repeat it
here with the extension, in the context of torsors, to all $N$.

The abstract theta group $\sG_L$ has an involution $\iota_L$
defined by $\iota_L(\lambda,P)=(\lambda, -P)$; this is the identity
on the central $\GG_m$. The concrete theta group
$\sG_{X,\Theta,N}$ also has such an involution, denoted by $\iota$
and defined by $\iota(f_P,P)=(f_P\circ[-1_A],-P)$.

\begin{definition} A level $(N,\vartheta)$ structure
$\psi:\sG_L\to\sG_{X,\Theta,N}$ is {\emph{symmetric}} if
$\psi\circ\iota_L=\iota\circ\psi$.
\end{definition}

Of course, if $N=2$ then every level $(N,\vartheta)$ structure
is symmetric.

For a symmetric level $(N,\vartheta)$ structure we can more nearly
pin down a normalization $(\alpha_P)$ of a Weil set $(f_P)$;
if $N$ is even, there is an ambiguity of $\pm 1$, and if $N$ is
odd there is no ambiguity.

As above, we can suppose that $\alpha_P\in\mu_N$ for all $P\in A[N]$.
The symmetry condition gives 
$\alpha_{-P}f_{-P}=\alpha_Pf_P\circ[-1_A],$
so that evaluation at one point $x_0$ of $X$ determines
the value of $\alpha_P/\alpha_{-P}.$ Since $\alpha_P\alpha_{-P}$
is determined by $\alpha_P\alpha_{-P}=\gamma(P,-P)$,
$\alpha_P$ is determined
uniquely if $N$ is odd and up to $\pm 1$ if $N$ is even.

\end{subsection}
\medskip
\begin{subsection}{Co-ordinates on moduli spaces}\label{co-ordinates}
We show now how to co-ordinatize first $A_{g,N,\vartheta}$ and then
$\sA_{g,N}$
over $\Z[\mu_N,1/N]$ by evaluating theta functions
(that is, sections of line bundles), 
without excluding characteristic $2$; of course,
this is very well known away from the prime $2$
and is due to Klein, Igusa and Mumford. To extend this
to include the prime $2$ is only slightly different, and we
do it by evaluating rational functions (which are ratios of theta functions)
on a torsor $X$ at the subscheme
$Fix_X$.

Start by identifying the vector space $H^0(X,\sO_X(N\Theta))$ with the
set of rational functions $g$ on $X$ whose polar divisor $(g)_\infty$ is at most 
$N\Theta$, together with $0$. There is a representation of $\sG_{X,N,\Theta}$ on 
$H^0(X,\sO_X(N\Theta))$,
given by $(f_P,P)(g)=f_P.t_P^*g$; since $\sG_{X,\Theta,N}$ 
is isomorphic to $\sG_L$, the representation theory of this latter group
tells us that this representation is irreducible.
Note that, if $g$ is a Weil function for $Q$, then
$(f_P,P)(g)$ is a Weil function for $P+Q$, so that the Weil functions 
lie in $H^0(X,\sO_X(N\Theta))$ and 
are permuted (up to scalar multiples) by $\sG_{X,\Theta,N}$. Therefore they span
$H^0(X,\sO_X(N\Theta))$.

A choice of level $(N,\vartheta)$ structure on $(X,\Theta)$
then determines a representation of $\sG_L$ on $H^0(X,\sO_X(N\Theta))$,
which makes $H^0(X,\sO_X(N\Theta))$ isomorphic to the unique irreducible 
representation $V_0$ of $\sG_L$, defined over $\Z[\mu_N,1/N]$, on which the 
central $\GG_m$ acts via homotheties. 

Let $V$ denote the $k$-vector space of functions $h:L\to k$.

\begin{lemma}
\part[i] $\sG_L$ acts on $V$ 
by $((\lambda,P)(h))(R)=\lambda.h(R-P)d(R,P)^{-1}.$

\part[ii] $L^\vee$ acts on $V$ by
$(\chi(h))(R)= \chi(R)h(R).$

\part[iii] There is a group structure $\tsG_L$ defined on
$\sG_L\times L^\vee$ by
$$((\lambda,P),\chi)((\mu,Q),\phi))=((\lambda,P)(\mu,Q)\chi(Q),\chi\phi).$$
Here the scalar $\chi(Q)$ is identified with an element of the
central $\GG_m$ in $\sG_L$.

\part[iv] Writing $((\lambda,P),\chi)=(\lambda,(P,\chi))$ exhibits
$\tsG_L$ as a central extension of $L\oplus L^\vee$
by $\GG_m$. As such $\tsG_L$ is 
isomorphic to $\sG_{L\oplus L^\vee}$.

\part[v] The actions of $\sG_L$ and $L^\vee$ on $V$
combine to give an action of $\G_{L\oplus L^\vee}$ that makes $V$ the
standard irreducible representation of $\G_{L\oplus L^\vee}$.
\begin{proof} Axiom checking.
\end{proof}
\end{lemma}

Suppose that $L,M$ are two symplectic $\Z/N$-modules with
Lagrangian decompositions
$L=L_1\oplus L_1^\vee$ and $M=M_1\oplus M_1^\vee$.
Then the space of maps $L_1\to k$  
is the standard irreducible representation $V_0$ of $\sG_L$;
let
$W_0$ denote that of $\sG_M$. Then $V_0\otimes W_0$ is naturally
isomorphic to the space of maps $L_1\oplus M_1\to k$, and so is the 
standard representation of $\sG_{L\oplus M}$. Then
$\sG_{L\oplus M}$ acts on the Segre embedding 
$\P(V_0)\times\P(W_0)\inj \P(V_0\otimes W_0)$.
The projections from $\P(V_0)\times\P(W_0)$ to
$\P(V_0)$ 
and to $\P(W_0)$ are $\sG_{L\oplus M}$-equivariant. Of course,
the central $\GG_m$ is acting trivially, so that $L\oplus M$
acts effectively on $\P(V_0)\times\P(W_0)$, while it acts
on $\P(V_0)$ via the projection $L\oplus\to L$
and on $\P(W_0)$ via $L\oplus\to M$. Also, the projection
$\P(V_0)\times\P(W_0)\to\P(V_0)$ factors through
$\P(V_0)\times\P(W_0)/M$, and
$\P(V_0)\times\P(W_0)\to\P(W_0)$ factors through
$\P(V_0)\times\P(W_0)/L$.

So, if we take $M=L^\vee$, then $V$ is isomorphic to $V_0\otimes W_0$
as representations of $\G_{L\oplus L^\vee}$.

Denote by
$\sX\to\sA_{g,N,\vartheta}$ the pullback of the universal principal
symmetric abelian torsor.

\begin{theorem} (Klein, Igusa, Mumford)
There is a morphism $\phi:\sX\to\P(V)$ defined by
$x\mapsto (\tf_P(x))_{P\in A[N]}$
whenever $(X,\Theta)$ is a principal symmetric abelian torsor 
and a level $(N,\vartheta)$
structure is chosen on $A=\Aut_X^0$. The $(\tf_P)_{P\in A[N]}$
are the normalized Weil functions determined by the level 
$(N,\vartheta)$ structure.
\noproof
\end{theorem}

This is nothing more than a very slight variant of the
use, which goes back to Klein, 
of rigidified global sections of line bundles over 
abelian varieties to construct moduli spaces for 
polarized abelian varieties $A$ with level structure. 
Note that we include characteristic $2$, but remember that $N$ is always prime
to the characteristic and $N\ge 2$.

\begin{theorem}\label{morphism to Hilb} Suppose that $N\ge 3$.

\part [i] 
There is a morphism
$\chi:\sA_{g,N,\vartheta}\to \Hilb_{\P(V),2^{2g}}$ given by taking
the image of $Fix_X$ in $\P(V)$ under $\phi$,
where $\Hilb_{\P(V),2^{2g}}$ is the Hilbert scheme
that parametrizes $0$-dimensional subschemes of degree $2^{2g}$ in $\P(V)$.

\part[ii] There is a commutative square
$$\xymatrix{
{\sA_{g,N,\vartheta}}\ar[r]^{\chi}\ar[d] & {\Hilb_{\P(V),2^{2g}}}\ar[d]\\
{\sA_{g,N}}\ar[r]^-{\phi} & {\Hilb_{\P(V),2^{2g}}/L^\vee.}
}$$

\part[iii] These morphisms are locally closed embeddings if $N\ge 4$.

\begin{proof} The existence of $\chi$ follows at once from
the preceding discussion. To see that it is a locally closed embedding
we appeal to 
the argument on pp. $199-200$
of \cite{MuTataIII} (which, as the authors there
state explicitly, works for all $N\ge 4$, whether even or odd).
This is to the effect
that, if we have a principal symmetric abelian torsor $(X,\Theta)$ over any base $S$ in which
$N$ is invertible, then for $N\ge 4$ the moduli point determines
the torsor, 
since $X$ can be recovered as the intersection of the quadrics through
the $\sG_L$-orbit of any given point of $Fix_X$.
\end{proof}
\end{theorem}


\begin{remark} 
\part[i] Going further, 
we can pass from the Hilbert scheme to the Chow scheme
$\Chow_{\P(V),2^{2g}}$ of $0$-cycles with degree $2^{2g}$,
which is just the symmetric product
$\P(V)^{(2^{2g})}$. That is,
there are morphisms $\sA_{g,N,\vartheta}\to \Chow_{\P(V),2^{2g}}$ 
and $\sA_{g,N}\to\Chow_{\P(V),2^{2g}}/L^\vee$,
and these morphisms are also locally closed embeddings
if $N\ge 4$. The homogeneous co-ordinates on this Chow scheme
are multi-symmetric functions of the co-ordinates on $\P(V)$,
so more explicit than those on the Hilbert scheme.

\part[ii] These morphisms to the Chow scheme can be defined
directly, without going to the Hilbert scheme; this approach
has the advantage of working also when $N=2$.

\part[iii]
In order to prove that these morphisms are locally closed embeddings
when $N=3$, we would need to know that there is an 
integer $d\ge 2$ such that $X$ can be recovered as the 
intersection of the $d$-ics through $Fix_X$. However, 
I do not know a proof of this.
\end{remark}

The variety that parametrizes the 
irreducible subrepresentations of $\sG_L$ in $V$, 
every one of which is isomorphic to $V_0$, 
is a copy of $\P(W_0)$.
In fact, $\P(V)$ contains a copy of the
Segre embedding of $\P(V_0)\times\P(W_0)$,
and $\P(V_0)\times\P(W_0)$ is the union of the projective spaces $\P(V_t)$
as $V_t$ runs over the family of irreducible subrepresentations of $\sG_L$ in $V$.

Since, for each geometric point $s$ of $\sA_{g,N,\vartheta}$,
the linear span of the image in $\P(V)$ of the torsor $\sX_s$
is $\P(V_t)$ for some $t$,
the image of $\sX$ under $\phi$ 
lies in $\P(V_0)\times\P(W_0)$. This gives the following assertion.

\begin{proposition} The morphism
$\chi:\sA_{g,N,\vartheta}\to\Hilb_{\P(V),2^{2g}}$ factors
(by abuse of notation) as a morphism
$\chi:\sA_{g,N,\vartheta}\to\Hilb_{\P(V_0)\times\P(W_0),2^{2g}}$.
\noproof
\end{proposition}

This leads to further quasi-projective realizations
of $\sA_{g,N,\vartheta}$ and $\sA_{g,N}$, as follows.

\begin{theorem}\label{morphism to projective space} 
\part[i] Suppose that $N\ge 3$. Then there 
are morphisms $\pi_1:\sA_{g,N,\vartheta}\to \P(W_0)$,
$\pi_2:\sA_{g,N}\to \Hilb_{\P(V_0),2^{2g}}$
and $\pi_3:\sA_{g,N}\to \Chow_{\P(V_0),2^{2g}}$.
\part[ii] $\pi_1$ and $\pi_3$ are also defined when $N=2$. 
\begin{proof}
Define
$\pi_1$ to the composite of $\phi$ with projection to $\P(W_0)$ 
(that is, $\pi_1$ takes a torsor $X$ 
to the linear span of its image under $\phi$).
Note that this argument works for all $N\ge 2$.

Suppose that $N\ge 3$.
To construct $\pi_2$ and $\pi_3$, note that
there is an open subscheme $\Hilb^0$ of $\Hilb_{\P(V),2^{2g}}$
such that the projection $pr_1:\P(V_0)\times\P(W_0)\to\P(V_0)$
induces a morphism $q_1:\Hilb^0\to \Hilb_{\P(V_0),2^{2g}}$.
Since, for every geometric point $s$ of $\sA_{g,N,\vartheta}$,
the image of the torsor $\sX_s$ in $\P(V)$ lies in one of
the copies of $\P(V_0)$ parametrized by $\P(W_0)$, the
image of 
$\chi:\sA_{g,N,\vartheta}\to\Hilb_{\P(V_0)\times\P(W_0),2^{2g}}$
lies in $\Hilb^0$. So composing $\chi$ with $q_1$ gives a morphism
$\tpi_2:\sA_{g,N,\vartheta}\to\Hilb_{\P(V_0),2^{2g}}$. Since
$L^\vee$ acts trivially on $\P(V_0)$, it follows that
$\tpi_2$ factors though a morphism
$\pi_2:\sA_{g,N,\vartheta}/L^\vee=\sA_{g,N}\to \Hilb_{\P(V_0),2^{2g}}$.

Finally, suppose that $N\ge 2$. Then the same argument, but working directly
with the Chow scheme, goes through for the construction of $\pi_3$.
\end{proof}
\end{theorem}

\begin{remark} The morphisms $\pi_2$ and $\pi_3$ can be summarized
as follows. Start with a principal symmetric abelian torsor $(X,\Theta)$ and a level
$N$ structure $\phi:L\to A[N]$,
corresponding to a geometric point $s$ of $\sA_{g,N}$.
Lift $s$ to a geometric point $\ts$ of $\sA_{g,N,\vartheta}$.
Then $\ts$ determines an isomorphism $H^0(X,\sO_X(N\Theta))\to V_0$.
The morphism $X\to\P(V_0)$ (an embedding if $N\ge 3$) restricts to
$Fix_X\to\P(V_0)$ and $\pi_2(s)$ or $\pi_3(s)$ is the point in 
$\Hilb_{\P(V_0),2^{2g}}$ or $\Chow_{\P(V_0),2^{2g}}$
corresponding to the image of $Fix_X$.
\end{remark}

It is not clear whether $\pi_1$ is an embedding,
while, for $N\ge 4$, $\pi_2$ is an embedding: the torsor $X$ can 
then be recovered as the intersection of the quadrics through $Fix_X$.

Instead of taking the image of $Fix_X$ in the Hilbert or Chow scheme
we can take the cycle-theoretic image of $X$ itself.

\begin{theorem} \part[i] Taking the cycle-theoretic image of $X$ in
$\P(V_0)\times\P(W_0)$ defines a morphism $\omega_1$ from
$\sA_{g,n,\vartheta}$ to the Chow scheme
$\Chow_{\P(V_0)\times\P(W_0),g,N^gg!}$
of cycles of dimension $g$ and degree $N^gg!$ in $\P(V_0)\times\P(W_0)$,
where ``degree'' refers to the degree in the ambient $\P(V)$.

\part[ii] Composing $\omega_1$ with the projection
$\P(V_0)\times\P(W_0)\to\P(V_0)$ defines a morphism
$\omega_2:\sA_{g,N}\to \Chow_{\P(V_0),g,N^gg!}$.
\noproof
\end{theorem}

This morphism $\omega_2$ is an embedding when $N\ge 3$, but is not quasi-finite when
$N=2$ and $g\ge 2$, because, if $X$ is, as a polarized torsor, a product $Y\times E$
with $E$ a curve of genus $1$, 
then this cycle cannot detect the isomorphism class of $E$.

Up to this point, in order
to remove the difference between level $(N,\vartheta)$ structures
and level $N$ structures, we have taken quotients by $L^\vee$.
Another way of doing this, at least when $N$ is odd, is to use
symmetric level $(N,\vartheta)$ structures,
as defined above, and the stack
$\sA_{g,N,\vartheta,symm}$ of ppav's with a symmetric
level $(N,\vartheta)$ structure; when $N$ is even this leads
to a level $(N,2N)$ structure. 
%
%

\begin{theorem} 
\part[i] There is a closed embedding
$\sA_{g,N,\vartheta,symm}\to\sA_{g,N,\vartheta}$.

\part[ii] The composite $\sA_{g,N,\vartheta,symm}\to\sA_{g,N}$ 
is an isomorphism if
$N$ is odd and is a quotient by the group $L^\vee/2L^\vee$ 
of quadratic characters of $L$ if $N$ is even.

\part[iii] If $N$ is even then $\sA_{g,N,\vartheta,symm}$
is isomorphic to the stack $\sA_{g,N,2N}$ of ppav's
with level $(N,2N)$ structure.
\begin{proof} This is an immediate consequence of \ref{theta group automorphism}
and the definition of level $(N,2N)$ structure.
\end{proof}
\end{theorem}

\begin{corollary} \part[i] If $N$ is odd, there are locally closed embeddings
$\sA_{g,N}\to\Hilb_{\P(V),2^{2g}}$ and $\sA_{g,N}\to\Chow_{\P(V),2^{2g}}$.

\part[ii] If $N$ is even, there are locally closed embeddings
$\sA_{g,N,2N}\to\Hilb_{\P(V),2^{2g}}$ and $\sA_{g,N,2N}\to\Chow_{\P(V),2^{2g}}$.
\noproof
\end{corollary}
\end{subsection}
\medskip
\begin{subsection}{Odd and even points (characteristic not $2$)}\label{char not 2}

As already mentioned, in order to be able to use the theta-constants
(that is, theta functions evaluated at a single point rather than
the entire subscheme $Fix_X$)
as homogeneous co-ordinates on the moduli space
we must exclude characteristic $2$.

Denote by $\sN_g=\sN_g^+\coprod\sN_g^-$ the stack of principally polarized
abelian $g$-folds $A$ with a symmetric divisor $\Theta$
that defines the polarization; $\sN_g^+$ or $\sN_g^-$
is the connected component
where $\mult_0\Theta$ is even or odd (\cite{FC}, p. 132). 
The forgetful map $\sN_g\to\sA_g$ is a torsor under $\sU[2]$
and $\sN_g^\pm\to\sA_g$ is finite and {\'e}tale of degree $2^{g-1}(2^g\pm 1)$.
We identify $\sN_g$ with the stack of principal symmetric abelian 
torsors $(X,\Theta)$ with a choice of point $x$ in $Fix_X$; then
$\sN_g^\pm$ is the open substack where the multiplicity of $x$
on $\Theta$ is even or odd.
\smallskip

Let $\sG_L^*=\sG_L\rtimes\langle\iota_L\rangle$ be the semi-direct product.

\begin{lemma} The standard $N^g$-dimensional representation $V_0$ of
$\sG_L$ extends uniquely to a representation of $\sG_L^*$.
\begin{proof} Choose a decomposition $L=L_1\oplus L_2$ into Lagrangian
submodules, so that $L_2$ is identified with $L_1^\vee$. Then $V_0$
has a basis $\{e_l\vert l\in L_1\}$ with respect to which $L_1$ acts by
$m(e_l)=e_{l-m}$, $L_2$ acts by $\chi(e_l)=\chi(l)e_l$ and the central
$\GG_m$ acts by $\lambda(e_l)=\lambda e_l$. We extend this by defining
$\iota_L(e_l)=e_{-l}$.
\end{proof}
\end{lemma}

As before, $V=V_0\otimes W_0$, where $W_0$ is the trivial $N^g$-dimensional
representation of $\sG_L^*$.

Suppose that $N$ is odd, and consider the stacks $\sN_{g,N}^\pm$.
The fixed point locus of $\iota_L$ on $\P(V)$ is the disjoint union
of two linear subspaces, $\P(V^+)$ and $\P(V^-)$, where $V^\pm$ is the
$\pm 1$-eigenspace of $\iota_L$.

\begin{theorem} There are locally closed embeddings
$\sN_{g,N}^\pm\to\P(V^\pm)$.
\begin{proof} These embeddings are defined by choosing a symmetric
level $(N,\vartheta)$ structure on the pointed principal symmetric abelian torsor $(X,\Theta,x)$
and then taking the point $(\tf_P(x))_{P\in A[N]}$ in $\P(V)$.
This point lies in $\P(V^\pm)$, as appropriate.
\end{proof}
\end{theorem}

Before considering further the case where $N$ is even, we recall
some results of Igusa \cite{I} and Mumford \cite{MuTh}.
Then $Fix$ is the
set $T$ of {\emph{theta characteristics}} and $A[2]$ is the group $P$ of
{\emph{period characteristics}}. We summarize some of the relevant linear algebra
over $\Z/2$;
Igusa (\cite{I}, p. 209 {\it{et seq.}}) 
uses multiplicative notation, while other sources use additive
notation. We shall use additive notation.

Write $\eta=\log e_2$,
where $e_2$ is the Weil pairing on $A[2]$ and
$\log:\mu_2\to\Z/2$, with inverse $\exp$, 
takes $1$ to $0$ and $-1$ to $1$.
Then $T$ is identified (\cite{I}, p. 209) with the set of functions
$L:T\to\Z/2$ such that $L(u+v)+L(u)+L(v)=\eta(u,v);$
such an $L$ is a quadratic form. Its Arf invariant
$\eta(L)\in\Z/2$ is given by
$\eta(L)=\log(2^{-g}\sum_{r\in P}\exp L(r))$ (\cite{I}, p. 210.)
Igusa's $e(L)$ is then $e(L)=\exp\eta(L)$.

If $L$ is a theta characteristic, regarded
as a point on the torsor $X$, then define
$e_*(L)\in\Z/2$ to be the multiplicity $\mult_L(\Theta)$ of the divisor $\Theta$
at the point $L$, taken modulo $2$. 
This equals the logarithm of the the scalar obtained by evaluating
the isomorphism $\phi$ that defines the symmetry of $\sO(\Theta)$
at the point $L$, and it is also the Arf invariant $\eta(L)$. 
Igusa writes $e(L)=\exp\eta(L)$. Moreover, it
coincides, when $A$ is the Jacobian of a curve $C$,
with Mumford's definition of $e_*(L)$
in \cite{MuTh}, by the Riemann singularity theorem.
If $a,b\in A[2]$, then \cite{MuTh}
$$e_*(L)+e_*(L+a)+e_*(L+b)+e_*(L+a+b)=\eta(a,b),$$
so that if we define $f_L:A[2]\to\Z/2$ by
$$f_L(a)=e_*(L+a)+e_*(L),$$
then $f_L$ is a theta characteristic in Igusa's sense.
That is, $\eta(a,b)=f_L(a)+f_L(b)+f_L(a+b)$,
so that writing $L$ instead of $f_L$ gives the two
ways of seeing $T$. 

The Arf invariant
divides $T=Fix$ as $T=T^+\coprod T^-$,
where $T^+$ is the set of even points (those with Arf invariant $0$)
and $T^-$ is the set of odd points (those with Arf invariant $1$).
This is preserved under monodromy, so
that the stack $\sA_g$ has, over $\Sp\Z[\frac{1}{2}]$,
two finite covers $\sA_g^{odd}$ and $\sA_g^{even}$, of degrees
$2^{g-1}({2^g-1})$ and $2^{g-1}(2^g+1)$, and $\sA_g^\pm=\sN_g^\pm$.

The next lemma, due to Igusa (Lemma $23$ of \cite{I}, p, $214$), serves
to identify $X$ with $A$ once a level $2$ structure is chosen.
A level structure is an isomorphism $\psi:P\to (\frac{1}{2}\Z/\Z)^{2g}$,
where $(\frac{1}{2}\Z/\Z)^{2g}$ has the standard symplectic form
$e(m_1,m_2;n_1,n_2)=\exp (2\pi i (m_1{}^tn_2 -n_1{}^tm_2))$; here
$m_1,...,n_2\in (\frac{1}{2}\Z/\Z)^g$ are row vectors and
${}^tx$ is the transpose of $x$. A choice also of
an element $\delta$ in $T$ then defines 
$\psi\circ\delta:T\to (\frac{1}{2}\Z/\Z)^{2g}$. 

In \cite{MuEqnsI} Mumford defines, when $X=A$ and the theta
characteristics are identified with the period characteristics, 
$$\log e_*^{\sO(\Theta)}(a)=\mult_a(\Theta)-\mult_0(\Theta)\pmod 2$$
for $a\in A[2]$. In this case, if $L,M$ are theta-characteristics with $M=L+a$,
then 
$$\log f_L(a)=e_*(L+a)-e_*(L)=\log(e_*^{\sO(\Theta)}(L+a)/e_*^{\sO(\Theta)}(L)).$$

\begin{lemma}\label{igusa} (Lemma $23$ of \cite{I}, p, $214$)
A choice of level $2$ structure on $A$ determines a unique
even theta characteristic $\delta$ such that when we write
$\psi\circ\delta(\alpha)=(m_1,m_2)$, we get
$e(\alpha)=\exp(2\pi im_1{}^tm_2)$.
\noproof
\end{lemma}

So a choice of level $2$ structure rigidifies $A[2]$ and determines an
isomorphism $X\to A$ that takes $\delta$ o $0_A$, 
$Fix_X$ to $A[2]$, and so rigidifies $Fix_X$.
In particular, if $N$ is even, then a level $(N,2N)$ structure determines a level
$2$ structure, and so rigidifies $Fix_X$.
Since giving a symmetric $(N,\vartheta)$ structure 
is equivalent to giving a level $(N,2N)$ structure, we recover
the following famous result. 

\begin{theorem}(Igusa, Mumford)\label{igusa, mumford}
If $N$ is even, there is a locally closed embedding
$\sA_{g,N,2N}\to\P(V)$ obtained by sending the principal symmetric abelian torsor $(X,\Theta)$
with a symmetric level $(N,\vartheta)$ structure to the point
$(\tf_P(\delta))_{P\in A[N]}$. Here the $\tf_P$ are the normalized
Weil functions determined by the symmetric level structure.
\noproof
\end{theorem}

Later on, in \ref{thomae}, we shall see, when the base field is $\C$,
how to write the normalized Weil functions
in terms of theta functions with characteristics, which will make it
clear that this theorem is exactly the classical co-ordinatization
of $\sA_{g,N,2N}$ by the $N$'th powers of the theta-nulls.
\end{subsection}
\end{section}
\bigskip
\begin{section}{Determinantal formulae for theta functions on $\Jac^{g-1}$
and $\Jac^0$}\label{construct}
\label{Determinant}
\medskip

In this section the characteristic $p$ is arbitrary
and $C$ is a curve of genus $g\ge 2$.

The Jacobian morphism $j:\sM_g\to\sA_g$ that sends a curve $C\to S$ of genus $g$
to its principally polarized
Jacobian $\Jac^0 C\to S$ factors as $\sM_g\to\sX_g\stackrel{\beta}{\to}\sA_g$,
by sending $C\to S$ to $(X=\Jac^{g-1}_{C/S},\Theta)$, where $\Theta$ is the divisor
consisting of effective classes;
in fact, the morphism is most easily constructed this way.
The symmetry of $(X,\Theta)$ is a consequence of Grothendieck-Serre duality.

Our aim is to extend a simplified version 
of Fay's addition theorem
(Proposition 2.16 of \cite{F}, p.29), as specialized
by Matone and Volpato (\cite{MaV}, Proposition 4.1).
As already mentioned, this goes back to Klein. In these formulae
there are constants that depend on the choice of basis
of the space of sections of some line bundle (typically pluricanonical
or sesquicanonical). However, the construction of the normalized Weil
functions has the effect of removing these unknown constants.

For every $n\ge 1$ there is an Abelian sum map 
$$\alpha=\alpha_n:C^n\to\Jac^nC,$$
constructed as a composite $C^n\to C^{(n)}\to\Jac^nC$.
Denote the $(i,j)$ diagonal in $C^n$ by $\Delta_{ij}$.

 
\begin{proposition}\label{formula}
Take the ground field to be algebraically closed;
then $\Theta$ is a well defined cohomology class on each
$\Jac^nC$
and the classes
$(2g-2)\alpha^*\Theta$
and
$(g-1+n)\sum_i pr_i^*K_C-(2g-2)\sum_{i<j}\Delta_{ij}$
on $C^n$ are cohomologous.

\begin{proof} It is convenient to prove this over the moduli space $M_g$,
so that we can exploit the fact the the image of the monodromy group
acting on the $2g$-dimensional vector space
$H^1(C,L)$, where the cohomology $H^*$ has coefficients
in a field $L$ of characteristic zero (for example,
$\ell$-adic cohomology) is Zariski dense in the symplectic group
$Sp(H^1(C,L))$.

Then $\alpha^*\Theta$ is a class in $H^2(C^n,L)$
and the class $\alpha^*\Theta$ is invariant under the symmetric 
group $\frak S_g$ and under $Sp(H^1(C,L))$.
So suppose that $V$ is a $2g$-dimensional vector space over $L$ with a non-degenerate
symplectic form $\phi$. The First Main Theorem for the invariant theory of $Sp(V,\phi)$,
together with consideration of the symmetric group $\frak S_n$,
asserts that, for any positive integer $n$,
$$(\bigwedge^2(V^{\oplus n}))^{Sp(V)\times\frak S_n}=  L.\sum_i\phi_i + L.\sum_{i<j}\psi_{ij},$$
where $\phi_i$ is just $\phi$ on the $i$th copy $V_i$ of $V$ and 
$\psi_{ij}$ corresponds to the given isomorphism $V_i\to V_j$.

Hence $\alpha^*\Theta$ is cohomologous to $a\sum_i pr_i^*K_C +b \sum_{i<j}\Delta_{ij}$,
for some $a,b\in L$. To determine $a$ and $b$ we intersect (separately) with two curves,
$\Gamma_1$ and $\Gamma_2$.

\smallskip
\noindent First, take $\Gamma_1$ to be the vertical curve (``a co-ordinate axis'') 
consisting of
the set $\{(P_1,\ldots,P_{n-1},x)\}$, where the $P_i$ are fixed and distinct.
Then $\Gamma_1.\alpha^*\Theta= g$, by Poincar\'e's formula, while
$\Gamma_1.\sum_i pr_i^*K_C = 2g-2$ and $\Gamma_1.\sum_{i<j}\Delta_{ij}=(n-1)$.
So $g=(2g-2)a+(n-1)b$.

\smallskip
\noindent Second, take $\Gamma_2$ to be the diagonal $\{(x,\ldots,x)\}$.
Then $\Gamma_2\to \Jac^n C$ factors as
$$\Gamma_2\to\Jac^1 C\stackrel{[n]}{\to}\Jac^n C,$$
where $[n]$ is multiplication by $n$. Since $[n]^*\Theta =n^2\Theta$,
in an obvious abuse of notation, Poincar\'e's 
formula shows that $\Gamma_2.\alpha^*\Theta=n^2g$.

On the other hand, note that $\sum pr_i^*K_C =K_{C^n}$ and
$\Gamma_2.\sum pr_i^*K_C =n(2g-2).$ Put $\Delta_{12}=\delta$,
a copy of $C^{n-1}$.
So $\Gamma_2.K_{C^n}=n(2g-2)$ and 
$\Gamma_2.K_{\delta}=(n-1)(2g-2)$. The adjunction formula
then shows that $\Gamma_2.\delta= -(2g-2)$. It follows that
$n^2g=a.n(2g-2)-b.(2g-2).n(n-1)/2.$

Clearing denominators now gives the result.
\end{proof}
\end{proposition}

Now suppose that $m$ is any integer $m\ge 2$ 
and $n=(2m-1)(g-1)$, and take $\Theta$ to be the effective divisor
on $\Jac^nC$ defined by the canonical isomorphism
$\Jac^{g-1}C\to\Jac^nC$.

\begin{corollary}\label{determinant} Suppose that $P$ is a point on $A$,
corresponding to the invertible sheaf $\sL_P=\sO_C(D_P)$
of degree $0$ on $C$,
and that $\{\sigma_1^P,\ldots,\sigma_n^P\}$ is a basis
of $H^0(C,\sO(mK_C+D_P)$. Then the zero locus 
on $C^n$ of $\det(\sigma_i^P(z_j))$ equals
$\alpha^{-1}(\Theta_P)+\sum \Delta_{ij}$.

\begin{proof} We have just shown
that $(\det(\sigma_i^P(z_j)))_0$ and
$\alpha^{-1}(\Theta_P)+\sum \Delta_{ij}$ are cohomologous.
On the other hand, any point in $\alpha^{-1}(\Theta_P)$ is of the
form $\{z_1,\ldots,z_n\}$ with $\sum_1^{(2m-2)(g-1)}z_i$
linearly equivalent to $(m-1)K_C+D_P$.
That is, there is a non-zero section $\tau$ of $H^0(C,\sO((m-1)K_C+D_P))$
vanishing at $\sum_1^{(2m-2)(g-1)}z_i$. 
Pick a basis $\omega_1,\ldots,\omega_g$ of
$H^0(C,\sO(K_C))$, and then
extend $\{\tau\omega_1,\ldots,\tau\omega_g\}$
to a basis $\{\sigma_i'\}$ of $H^0(C,\sO(mK_C+D_P))$.
Then the $n\times n$ matrix
$(\sigma_i'(z_j))$ has a $g\times (2m-2)(g-1)$ block that
is identically zero. This forces $\det(\sigma_i'(z_j)) =0$,
and then the same is true of $\det(\sigma_i^P(z_j))$.
That is, $\alpha^{-1}(\Theta_P)+\sum \Delta_{ij}$ 
is contained in the member $(\det(\sigma_i^P(z_j))_0$ of 
the linear system $\vert \sum pr_i^*(mK_C+D_P)\vert$;
since they are cohomologous over the algebraic closure of $k$, they are equal.
\end{proof}
\end{corollary}

\begin{remark} This is a version, valid over any field, of 
Proposition 4.1 of \cite{MaV}. There the ground field is assumed to
be $\C$, but they also give explicit transcendental formulae for the
constants of proportionality that we do not.
\end{remark}

Now suppose that $n=2m(g-1)$, with $m\ge 1$.
Then $\Jac^n_C$ is canonically isomorphic to $\Jac^0_C$.
There is no theta divisor, but only its cohomology class
(over an algebraic closure of the ground field);
there is a canonically defined line bundle $\sL=\sL^\Delta(\lambda)$,
where $\lambda$ is the principal 
polarization on $\Jac^0_C$, that defines twice the theta
cohomology class $[\theta]$.

Suppose that $M, M'$ are
divisor classes
of degree $g-1$ on $C$ and that
$\{\sigma_1^M,\ldots,\sigma_n^M\}$ is a basis
of $H^0(C,\sO(mK_C+M))$. 
Denote by $\Theta_{-M}$ the image of $\Theta$ under the 
subtraction isomorphism
$\Jac^{g-1}_C\to\Jac^0_C$ defined by $M$.

\begin{proposition}
\part[i]\label{theta_char_det} 
The zero locus on $C^n$ of $\det(\sigma_i^M(z_j))$
is $\alpha^{-1}(\Theta_{-M})+\sum_{i<j}\Delta_{ij}$.

\part[ii] $\alpha^{-1}(\Theta_{-M})$ is linearly equivalent to
$\sum pr_i^*(mK_C+M)-\sum\Delta_{ij}$.

\part[iii] If $M$ is a theta characteristic, then
$2\Theta_{-M}$ is in the linear system $\vert \sL\vert$.

\begin{proof} The proof of \DHrefpart{i} is 
the same as that of \ref{determinant}, and
\DHrefpart{ii} follows at once.

Certainly $2\Theta_{-M}$ and $\sL$ define the same cohomology
class, and $\alpha^*(2\Theta_{-M})$ is linearly equivalent
to $\sum pr_i^*((2m+1)K_C)-2\sum\Delta_{ij}$. But on the generic
Jacobian, the only rational point is $0$, so that
$2\Theta_M$ and $\sL$ are linearly equivalent generically,
so everywhere.
\end{proof}
\end{proposition}

\begin{corollary}\label{cassels_flynn} 
If $n=2g-2$, then 
$\alpha^*\sL$ is linearly equivalent to
$3\sum_1^{2g-2}pr_i^*K_C-2\sum\Delta_{ij}$.
\noproof
\end{corollary}
\end{section}
\bigskip
\begin{section}{The group law on $A$ via chords and tangents
}\label{group law}\label{chord and tangent}
\medskip
Here we give a projective description of the group law on $A$,
analogous to the chord and tangent construction
for plane cubics.

Suppose that $(X,\Theta)$ is any principal torsor and that
$(A,\lambda)$ is the associated ppav. 
Recall that for $N\ge 3$, the linear system $\vert N\Theta\vert$
on $X$ is very ample and embeds $X$ in $\P^{N^g-1}$.

\begin{proposition} There is a description of the group law on $A$ 
in terms of a chord and tangent construction on the $\vert 3\Theta\vert$
image of $X$,
as described in the course of the proof.

\begin{proof} Given $P$ and $Q$ in $A$,
there is, because $H^0(X,\sO(\Theta_x))=1$ for all $x$ in $A$,
a unique hyperplane $H$ in $\P^{3^g-1}$ that contains $\Theta_P$ and $\Theta_Q$,
so that $H.X=\Theta_P +\Theta_Q +E$ for some divisor $E$. Then $E=\Theta_{-P-Q}$. 
Similarly, there is a unique hyperplane $H'$ containing $E$ and
$\Theta_0$, and $\Theta_{P+Q}$ is recovered as the residual
intersection of $H'$ with $X$.
\end{proof}
\end{proposition}

It is important to notice that this can also be interpreted as the statement that
there is a morphism $\pi$ from $A$ to the Grassmannian $G(2^g,3^g)$
such that each $P\in A$ is mapped to the $2^g$-dimensional space of linear forms
that vanish along $\Theta_P$ and the group law on $A$ has a projective description.
Obviously, $\pi$ separates points.

Recall that the principal polarization $\lambda$ on $A$ defines a totally
symmetric ample sheaf $\sL^\Delta_\lambda$ on $A$ that returns $2\lambda$.
Denote by $pr_1,pr_2:A\times A\to A$ the projections and
$\sigma:A\times A\to A$ the group law. We also let $\Theta$
denote the cohomology class on $A$ defined by $\Theta$, or $\lambda$.

\begin{proposition} The composite morphism $j:A\to\P^N$ defined by 
$\pi$ and the Pl{\"u}cker embedding of $G(2^g,3^g)$ in $\P^N$, where
$N={3^g\choose 2^g}-1$, separates points and is defined by a subsystem of
the complete 
linear system $\vert m\sL^\Delta_\lambda\vert$, where
$$m=-(pr_1^*\Theta)^{g-1}(3pr_2^*\Theta-\sigma^*\Theta)^{g+1}/2(g+1)!g!$$

\begin{proof} It is clear that $j$ separates points.
Put $V_3=H^0(X,\sO_X(3\Theta))$. We have
two morphisms $pr_2,\sigma:A\times X\to X$, where $\s$ is the action, and 
a morphism $pr_1:A\times X\to A$.

Put $\sM=pr_2^*\sO_X(3\Theta)\otimes\sigma^*\sO_X(-\Theta).$
This is naturally a subsheaf of $pr_2^*\sO_X(3\Theta)$; taking sections
gives a short exact sequence
$$0\to pr_{1*}\sM\to\sO_A\otimes V_3\to \sE\to 0$$
of locally free sheaves on $A$, where 
$\sE=pr_2^*\sO_X(3\Theta)/\sM$, by definition.

By construction, $\pi:A\to G(2^g,3^g)$ is defined by the quotient
$\sO_A\otimes V_3\to \sE\to 0$. 

Define $pr_{1*}\sM=\sF$. Then,
by Grothendieck-Riemann-Roch, $ch(\sF)=pr_{1*}(ch(\sM))$.
Since $ch(\sM)=\exp(c_1(\sM)$, we get
$c_1(\sF)=pr_{1*}(3pr_{2*}\Theta-\sigma^*\Theta)^{g+1})/(g+1)!$
On the other hand, we know, for example by considering the generic
principally polarized abelian variety, that
$\det\sE$ is isomorphic to $(\sL^\Delta_\lambda)^{\otimes m}$
for some $m$. To compute $m$, it is enough to determine
$\Theta^{g-1}.c_1(\sF)$, since $c_1(\sF)=-c_1(\sE)$
and $\Theta^{g-1}.c_1(\sE)=2g!m$.
This gives exactly the stated formula for $m$.
\end{proof}
\end{proposition}

The morphisms $pr_1,pr_2,\sigma:A\times A\to A$ can be made symmetric
by identifying $A\times A$
with the subvariety $\{(P,Q,R)\vert P,Q,R\in A, P+Q+R=0\}$
and replacing $\sigma$ by its composite with $[-1_A]$.

For example, if $g=2$, then $m=3$. Then $A$ is embedded in
$\P^{125}$ via a $6\Theta$ system, so is the intersection
of $90$ hyperplanes and $29.18=522$ quadrics.

To elaborate slightly,
suppose that $T$ is the homogeneous co-ordinate ring of $\P^{3^g-1}$, that
$S=T/J$ is the homogeneous co-ordinate ring of $X$ with respect to
$3\Theta$ and that $I_P,I_Q$ are the homogeneous prime ideals in $S$
of $\Theta_P, \Theta_Q$. There is a linear form $f$, unique up to scalars,
defining the linear span $\langle \Theta_P,\Theta_Q\rangle$; we need to
compute the vector space $W_R$ of 
linear forms $h$ that contain the divisor $\Theta_R$ defined by
$\Theta_P+\Theta_Q+\Theta_R=X\cap (f=0)$ (that is, $R=-P-Q$). 

Assume for the moment that $P\ne Q$. 
Then the primary decomposition of the ideal $(f)$ of $S$ is
$(f)=I_P\cap I_Q\cap I_R$, so the condition on the linear forms $h$ is that
$h.p$ should lie in $(f)$ for every $p\in I_P\cap I_Q$.
Gr{\"o}bner basis techniques \cite {BW} now permit a translation of this into
a polynomial formula for the Pl{\"u}cker co-ordinates
of $W_R$ in terms of the Pl{\"u}cker co-ordinates of $W_P$ and $W_Q$. 

This is of limited use without the ability to compute the vector spaces $W_P$ and
$W_Q$. For Jacobians, we can use the determinantal formula to make these computations.

Assume, for simplicity, that $C$ is non-hyperelliptic.
Suppose that $X=\Jac^{g-1}_C$ and $\Theta$ its theta divisor, as before.
We know that $\alpha^*\Theta$ is linearly equivalent to
$2\sum_i pr_i^*K_C-\sum_{i<j}\Delta_{ij}$.
Put $\Delta=\sum_{i<j}\Delta_{ij}$ and regard it as a closed subscheme
of $(\P^{g-1})^{3g-3}$ via the embeddings
$\Delta\inj C^{3g-3}\inj (\P^{g-1})^{3g-3}$. Take the subscheme
$Z_N=N\Delta$ of $C^{3g-3}$, so that the ideal sheaf $\sI_{Z_N}$
of $Z_N$ in $(\P^{g-1})^{3g-3}$ is $\sI_{Z_N}=\sI_\Delta^N+\sI_{C^{3g-3}}$.

Assume that either $p=0$ or $p>3g-3$, so that $(3g-3)!$ is invertible in $k$
and 
the symmetric group $\frak S=\frak S_{3g-3}$ is linearly reductive.
Let us say that a form in $3g-3$ variables is
\emph{semi-invariant of order $N$} if it transforms as
the $N$th power of the signature of $\frak S$. 
We shall write $V^{\frak S,N}$ for the space of order
$N$ semi-invariants in $V$.

\begin{proposition}\label{polynomials}\label{thetas on X} 
\part[i] $H^0(X,\sO_X(N\Theta))$ is naturally isomorphic to
the vector space $H^0(C^{3g-3},\sO(2N\sum pr_i^*K_C-N\sum\Delta_{ij}))^{\frak S,N}$
of order $N$ semi-invariants.

\part [ii] There is a natural inclusion from
$$(H^0((\P^{g-1})^{3g-3}, \sI_{Z_N}.\sO(2N,\ldots,2N))/
H^0((\P^{g-1})^{3g-3}, \sI_{C^{3g-3}}.\sO(2N,\ldots,2N)))^{\frak S,N}$$
to $H^0(X,\sO_X(N\Theta))$. This is an isomorphism if
either $p=0$ or $p>(3g-3)!$.

\begin{proof} 
Noether's theorem on canonical curves is that
the natural map $$H^0(\P^{g-1},\sO(r))\to H^0(C,\sO(rK_C))$$
is surjective for all $r\ge 1$. Then the same holds
for 
$$H^0((\P^{g-1})^{3g-3}),\sO(r,\ldots,r))\to H^0(C^{3g-3},\sO(r\sum_{i=1}^{3g-3} pr_i^* K_C)).$$
Taking semi-invariants preserves surjectivity
if $\frak S_{3g-3}$ is linearly reductive.

The proposition now follows
from the determinantal formula \ref{determinant} 
for the pull-back of the theta divisor
under the addition map $\alpha:C^{3g-3}\to\Jac^{g-1}_C$,
except for the statement about the order of the semi-invariants.
For this, it is enough to take $N=1$ and note the skew-symmetry
of the determinant $\det(\sigma_i(P_j))$, where $\{\sigma_i\}$
is a basis of $H^0(C,\sO(2K_C))$.
\end{proof}
\end{proposition}

So, given the defining equations for $C$, we can take $N=3$
and then have an explicit description
of $H^0(X,\sO_X(3\Theta))$. In terms of an arbitrary basis of this,
we can then
write down the morphism $\phi:X\to \P^{3^g-1}$ defined by $\vert 3\Theta\vert$.
That is, we can compute $\phi(P)$ for any point $P\in C^{3g-3}$ not in $\Delta$.
Since $X$ in $\P^{3^g-1}$ is the intersection of a known number of cubics,
we can compute $\phi(P)$ for many points $P$ and then find the cubic equations
that define $X$ by interpolation. Then, given $P\in A=\Jac^0_C$ and $x\in C^{3g-3}$
that represents a point of $\Theta$, it is possible to use the projective
geometry of $C$ to compute representatives of $P+x$ in $C^{3g-3}$.
Then we can compute $\phi(P+x)$ for many points $x\in \Theta$, and so
compute the vector space of hyperplanes through $\phi(\Theta_P)$. 
After that, the results above for an arbitrary principally polarized abelian
variety can be applied.

In particular, we can write down polynomial equations that determine
whether a given point $P$ on $A$ is $N$-torsion,
in terms of the Pl{\"u}cker co-ordinates of the space $W_P$.
We regard this as a construction of the $N$-torsion divisor classes on $C$.

Similarly, \ref{cassels_flynn} can be restated to give
an explicit description of the spaces $H^0(A,\sL^{\otimes r})$
of theta functions of order $2r$ on $A=\Jac^0_C$. Recall that
$\sL=\sL^\Delta(\lambda)$ is the natural totally symmetric line bundle on $A$
that lies in the $2\theta$ class.

\begin{proposition}\label{thetas on A} 
Put $\frak S=\frak S_{2g-2}$ and let $r\ge 1$ be an integer. 

\part[i] There is a natural isomorphism 
$$H^0(A,\sL^{\otimes r})\to H^0(C^{2g-2},\sO(3r\sum pr_i^*K_C-2r\sum\Delta_{ij}))^{\frak S}.$$

\part [ii] There is a natural inclusion from
$$(H^0((\P^{g-1})^{2g-2}, \sI_{Z_{2r}}.\sO(3r,\ldots,3r))^{\frak S}/
H^0((\P^{g-1})^{2g-2}, \sI_{C^{2g-2}}.\sO(3r,\ldots,3r)))^{\frak S}$$
to $H^0(A,\sL^{\otimes r})$.
This is an isomorphism if
either $p=0$ or $p>(2g-2)!$.

\begin{proof} The only thing to be checked is that $H^0(A,\sL^{\otimes r})$
consists of invariants rather than semi-invariants. For this, take
a theta characteristic $M$ 
and consider the determinant $\det(\sigma^M_i(z_j))$ of \ref{theta_char_det}.
It is semi-invariant, so its square is invariant and is
a section of $\alpha^*\sL$.
\end{proof}
\end{proposition}
\begin{remark} This description of vector spaces of theta functions on Jacobians
(either $\Jac^{g-1}$ or $\Jac^0$) relies upon already having constructed the 
Jacobian; for example, without that, it would not be clear that the various
linear systems on $C^{(3g-3)}$ or $C^{(2g-2)}$ had no unassigned
base points. That is, it does not lead to another elementary projective
construction of the Jacobian such as that given by Anderson \cite{An}.
\end{remark}
\end{section}
\bigskip
\begin{section}{Computing the $N$-torsion}\label{computing torsion}
\medskip
As already mentioned, the chord and tangent construction
of the group law on $A=\Jac^0_C$ described above
specializes to a description of the
multiplication by $N$ map $[N_A]$, and therefore the kernel $A[N]$.
However, a more direct calculation of $[N_A]$ and the inversion $[-1_A]$
is also possible. The idea is to use the explicit descriptions
of spaces of theta functions given by \ref{thetas on X} and \ref{thetas on A}
and then derive the equations of $A$ and $X$ and formulae for $[N_A]$
by {\emph{interpolation}}.

Interpolation is the process that, given a projective
$k$-variety $V$ and a morphism $\phi:V\to\P^n$ such that, first, $\phi(x)$
can be calculated for all fields $K$ containing $k$ and for all
$K$-points $x$ of $V$ and, second, for every
integer $d\ge 1$ the dimension of the $k$-vector space $H^0(\P^n,\sI_{\phi(V)}(d))$
of degree $d$ hypersurfaces that contain the image $\phi(V)$ is known,
produces a $k$-basis of that space. So, if also the ideal
$\sI_{\phi(V)}$ is known to be generated by quadrics, then interpolation
produces, in a finite calculation, the equations that define $\phi(V)$.
The most efficient way of describing it is to take $x$ to be
defined over some finite extension of the function field of $V$
(for example, the generic
point of $V$) and then find the degree $d$ hypersurfaces defined over $k$
that contain $\phi(x)$.

Put $\sL=\sL^\Delta_\lambda$, the natural $2\Theta$ line bundle on $A$
that is determined by the principal polarization $\lambda$.
Put $4^g-1=r$; then,
according to \ref{thetas on A} there is a basis
$z_0,...,z_r$ of $H^0(A,\sL^{\otimes 2})$
whose pull back to $C^{2g-2}$ can be computed explicitly. Fix $N\ge 2$
and generic points $\eta_1,...,\eta_{2g-2}$ on $C$.
Take $x=(\eta_1,...,\eta_{2g-2})$, a generic point of $C^{2g-2}$.
Put $L=k(C^{2g-2})$.

%

Let $\psi:C^{2g-2}\to A\to\P^r$
be the composite of the abelian sum $\alpha$ and the embedding $A\to\P^r$
defined by $H^0(A,\sL^{\otimes 2})$.
\medskip

\noindent (1) Since in this embedding
$A$ is known to be an intersection of quadrics and we know that
the quadrics in $\P^r$ cut out a complete linear system on $A$,
it follows that $\dim H^0(\P^r,\sI_A(2))={{r+2}\choose{2}}-8^g$.
Hence 
we can interpolate the quadrics through the point $p=\psi(x)$
to write down (a $k$-basis of the vector space of) the quadrics through $A$.

\noindent (2) It is known that $H^0(\P^r,\sI_A(N^2))$ is generated by
$H^0(\P^r,\sI_A(2))$. So we can write down a $k$-basis of $H^0(\P^r,\sI_A(N^2))$
and then a $k$-basis $w_0,...,w_s$ of a complement $V$ to $H^0(\P^r,\sI_A(N^2))$
inside $H^0(\P^r,\sO(N^2))$. Of course, $s=(4N^2)^g-1$ and $V$
represents $H^0(A,\sL^{\otimes 2N^2})$.

\noindent (3) 
Compute a point $y\in C^{2g-2}$ defined over a finite extension $L'$ of $L$
such that $y+(N-1)K_C$ is linearly equivalent to $Nx$ as divisors on $C\otimes L'$, 
so that
if $p=\alpha(x)$ and $q=\alpha(y)$, then $q=Np$ in $A$.

\noindent (4)
Since $[N_A]^*\sL$ is isomorphic to $\sL^{\otimes N^2}$, the map $[N_A]:A\to A$
is defined by formulae $z_j(Np)=f_j(z_0(p),...,z_r(p))$,
where the $f_j$ are linear combinations of $w_0,...,w_s$
with undetermined coefficients.

\noindent (5) Then $z_j(q)=f_j(z_0(p),...,z_r(p))$ for every 
$j$. Since we can calculate $z_j(p)$ and $z_j(q)$, we can determine
the $f_j$ as a linear combination of $w_0,...,w_s$.
That is, we now have formulae for the map
$[N_A]$.

\noindent (6) Take a $k$-point $x_0\in C^{2g-2}$ that is linearly equivalent
to $K_C$ and compute its image in $A$: this is the origin $0_A$.
By combining the quadratic equations for $A$
with the equation $Np=0_A$, we have equations for the $N$-torsion
group $A[N]$.

The $2$-torsion can also be computed as the fixed locus in $A$ (or $X$)
under $[-1_A]$; this is a little quicker.

Given $x_i\in C^{2g-2}$, compute $y_i\in C^{2g-2}$ with $y_i$
linearly equivalent to $2K_C-x_i$. Then compute the images in $\P^r$
of the points $x_i$ and $y_i$; they span a line that is preserved under
the involution $[-1_A]$ of $\P^r$. After sufficiently many of these
lines have been computed, we can determine $[-1_A]$ as a linear
transformation of $\P^r$, and so determine the two linear subspaces
whose union is the fixed locus in $\P^r$. This gives equations
for the inverse image in $C^{2g-2}$ of $A[2]$.
\medskip

It is also worth remembering that, when 
$g=3$ and $g=4$,
Coble showed \cite{C} how
to write down non-hyperelliptic curves together with the $2$-torsion on their Jacobians
(or, equivalently, their odd theta-characteristics),
provided that when $g=4$ the curve has a vanishing even theta-null;
this happens exactly for non-hyperelliptic
genus $4$ curves whose canonical
model lies on a quadric cone. The curves are constructed
as the ramification divisor of a double cover $S\to\P^2$, where
$S$ is a del Pezzo surface of degree $2$ or $1$ accordingly.
This should be compared
to Weber's determination of the bitangents of a plane quartic
as used by Ritzenthaler \cite{R}.

Coble's results are as follows.
Recall that a set $Z$ of $r$ points in $\P^2$ is {\emph{in general position}}
if they are distinct, no $3$ are collinear, no $6$ on a conic and no $8$
on a cubic that is singular at one of the $8$.
Suppose that $r\le 8$; then blowing up $\P^2$ at $Z$ gives a del Pezzo
surface $S$ of degree $9-r$. The complete anti-canonical system $\vert{-}K_S\vert$
is the system $\vert 3H-Z\vert$ of cubics in $\P^2$ through $Z$.
\smallskip

Suppose that $r=7$. 
Then $\vert{-}K_S\vert$ defines a finite double cover $\pi:S\to\P^2$ whose branch
locus $B\subset\P^2$ is a smooth quartic, and every smooth quartic
arises in this way. The $28$ bitangents are the images of the $56$
exceptional curves of the first kind (``lines'') on $S$.

Regard $\pi$ as the rational map $\pi:\P^2-\to\P^2$ defined by the linear
system $L=\vert 3H-Z\vert$. Then we get a birational model of a genus $3$
curve $C$ in the source
$\P^2$ as the ramification curve (that is, the locus where the gradient vectors
of the cubics in $L$ span a space of rank $1$); $C$ is birational
to a member $\barC$
of $\vert 6H-2Z\vert$. The $7$ exceptional curves $E_i$ in $S$ over the points in $Z$
and the $21$ lines $L_{jk}$ in $\P^2$ spanned by pairs of points in $Z$ represent
the bitangents, so we get the $28$ odd theta-characteristics $P_i+Q_i$ as follows:
if $i=1,...,7$, then $P_i+Q_i$ is represented by a node on $\barC$, and otherwise
$P_{jk}+Q_{jk}$ is the residual intersection of $L_{jk}$ with $\barC$
(\cite{C}, pp. 158-9).
If $F_1,F_2,F_3$ span $\vert 3H-Z\vert$ and $X_1,X_2,X_3$ are homogeneous co-ordinates
on $\P^2$, then the equation of $\barC$ is $\det(\partial F_i/\partial X_j)=0$.
\smallskip

If $r=8$, then $\vert{-}2K_S\vert$ defines
a finite double cover $\pi:S\to Q$, where $Q$ is a quadric cone in $\P^3$.
In this case the branch locus is a smooth genus $4$ curve $C\subset Q$,
the intersection of $Q$ with a cubic surface, together with the vertex of $Q$.
Every non-hyperelliptic curve of genus $4$ with a vanishing even theta-null
arises in this way. Again, regard
$\pi$ as the rational map $\pi:\P^2-\to Q\subset\P^3$ defined by the linear
system $L=\vert 3H-Z\vert$. This time, $C$ is birational
to a member $\barC$
of $\vert 9H-3Z\vert$ in the source $\P^2$.
If $F_1,F_2$ span $\vert 3H-Z\vert$ and $F_1^2,F_2^2,F_3$ span
$\vert 6H-2Z\vert$, then again the equation of $\barC$ is
$\det(\partial F_i/\partial X_j)=0$.
The $120$ odd theta-characteristics arise as follows. There are $8,28,56,28$
that correspond, respectively, to the points of $Z$ (the triple points of $\barC$),
the lines spanned by pairs of points in $Z$, the conics through $5$-tuples in $Z$
and one half of the $56$ cubics through $6$-tuples in $Z$ and singular at one
more point in $Z$. This half should be chosen to give a complete set of
representatives under the Bertini involution $\iota$ of $S$, which is the
Cremona transformation of $\P^2$ defined by the linear system $\vert 17H-6Z\vert$
(\cite{C}, p. 209).
\end{section}
\bigskip
\begin{section}{Explicit Weil functions on $\Jac^{g-1}$}\label{weil}
\medskip

Suppose that $X=\Jac^{g-1}_C$ and $N$ is prime to $p$; then, 
given an $N$-torsion point $P$ on $A$, which we regard
as an invertible sheaf of order $N$ on $C$, 
we can construct Weil functions when $X=\Jac^{g-1}_C$,
or their pull-backs to $C^{3g-3}$, as follows. 

Start by constructing a basis $\{\sigma_i^P\}$ of $H^0(C,\sO(2K_C)\otimes P)$
for every $P\in A[N]$, where each such point is regarded as a line bundle on $C$.

\begin{proposition}
The rational function $f_P$ on $C^{3g-3}$ defined by
$$f_P(z)=\det(\sigma_i^P(z_j))^N/\det(\sigma_i(z_j))^N$$
is the pull-back of a rational function on $X$ which is a Weil function there
for the point $P$ in $A[N]$.

\begin{proof} Since the given right-hand expression is ${\frak S}_{3g-3}$-invariant,
and the addition morphism $\alpha:C^{(3g-3)}\to\Jac^{3g-3}C$ satisfies
$$H^0(C^{(3g-3)},\alpha^*\sM)=H^0(\Jac^{3g-3}C, \sM)$$ 
for all line bundles
$\sM$ on $\Jac^{3g-3}C$, the result follows from the previous discussion.
\end{proof}
\end{proposition}

\begin{remark} We could use $\sO(nK_C)$ for any $n\ge 2$, as do Matone and Volpato.
For most of our purposes,
replacing determinants of order $3g-3$ by ones of order $(2n-1)(g-1)$
doesn't seem to help, although in section \ref{spinor} we do consider
such higher-order forms.
\end{remark}
\end{section}
\bigskip

\begin{section}{Comparison with theta functions}\label{theta}
\medskip
It is easy to translate this into the 
complex analytic language of theta functions
with characteristics. 

So suppose that 
$(X,\Theta)$ and $(A,\lambda)$ are as usual, but over $\C$.
Then the groups $H_1(X,\Z)$ and $H_1(A,\Z)$ are identified,
with their symplectic pairings, and the vector spaces
$H^0(X,\Omega^1_X)$ and $H^0(A,\Omega^1_A)$
are also identified.

A choice of symplectic
basis of $H_1(X,\Z)$, from which a normalized basis of
$H^0(C_\C,\Omega^1)$ and period matrix $\tau\in\frak H_g$, Siegel's upper half-space, 
are constructed, leads to identifications $H^0(A,\Omega^1_A)=\C^g$
and $A=A_\tau=\C^g/\Lambda_\tau$, where
$\Lambda_\tau =\tau \Z^g\oplus\Z^g$.
We then identify
$H^1(X,\Z/N)$, $H^0(A_\tau,\Z/N)$ and $A_\tau[N]$ with $\frac{1}{N}\Lambda/\Lambda$;
the cup product on $H^1(X,\Z/N)$, the symplectic form on $H^1(A_\tau,\Z/N)$
given by the polarization on $A_\tau$ and the Weil pairing on $A_\tau[N]$
are then identified with
$$\langle\ ,\ \rangle:\frac{1}{N}\Lambda/\Lambda\times \frac{1}{N}\Lambda/\Lambda\to\mu_N$$
given by $\langle P,Q\rangle = d(P,Q)/d(Q,P),$ where 
$d(\tau a+b,\tau e +f)=\exp(2\pi i N\ {}^teb)$
and $\zeta=\exp(2\pi i/N)$.

The chosen symplectic basis of $H_1(X,\Z)$ determines
a level $2$ structure on $(A,\lambda)$, or on $(X,\Theta)$.
According to \ref{igusa}, a level $2$ structure determines an even theta-characteristic
$\delta$ with certain properties.
Recall that for any real vectors $a,b$, the theta function
with characteristics $a,b$ is defined by
$$\theta{a\brack b}(z,\tau)=\sum_{n\in\Z^g}\exp(\pi i(2{}^t(n+a)(z+b)+{}^t(n+a)\tau(n+a)).$$
All vectors are column vectors of size $g$, ${}^tx$ is the transpose of $x$
and $z\in\C^g$.

Then Riemann's theorem can be stated as follows.

\begin{theorem}\label{Riemann} The isomorphism $X\to A$ determined
by $\delta$ takes $\Theta$ to the zero locus of the theta function
$\theta{0\brack 0}(z,\tau)$.\noproof
\end{theorem}

Recall that $\theta{a\brack b}(z,\tau)$ is quasi-periodic, with factors of automorphy
$$\theta{a\brack b}(z+p,\tau)=\exp(2\pi i\ {}^t a p)\theta[a,b](z,\tau),$$
$$\theta{a\brack b}(z+\tau q,\tau)= \exp(-2\pi i\ {}^t b q-\pi i\ {}^tq\tau q-2\pi i\ {}^tq z).
\theta{a\brack b}(z,\tau).$$
Moreover, 
$$\theta{a\brack b}(z,\tau)=\exp(\pi i\ {}^ta\tau a+2\pi i\ {}^ta(z+b))
\theta{0\brack 0}(z+\tau a+b,\tau),$$
so that $(\theta{a\brack b})_0=t^*_{-\tau a-b}(\theta{0\brack 0})_0.$

More generally, if $N$ is even, then a level $N$ structure on $A$ determines
$\delta$. So suppose $N$ even and fix, in fact, a level $(N,2N)$ structure on $A$.

Define $g{a\brack b}=\theta{a\brack b}/\theta{0\brack 0}$. 
Assume that $a,b\in(\frac{1}{N}\Z)^{g}$;
then $\tau a +b$ represents an $N$-torsion point on $A=A_\tau=\C^g/(\Z^g+\tau\Z^g)$,
and $g{a\brack b}^N$ is a meromorphic function on $A_\tau$ with
divisor $(g{a\brack b}^N)=N\Theta_{-\tau a -b}-N\Theta_0$.
Moreover, $\theta{a+A\brack b+B}/\theta{a\brack b}$ is an $N$th root of unity
if $A,B\in\Z^g$,
so that $\theta{a\brack b}^N$ 
depends only on the classes of $a,b$ in $(\frac{1}{N}\Z)^{g}/\Z^g$.

\begin{lemma} The set $\{g{-a\brack -b}^N\}_{\tau a+b\in A[N]}$
is a normal Weil set.
\begin{proof} Define $\phi{a\brack b}=g{-a\brack -b}$. Then 
$$\phi{e\brack f}.t_{\tau e + f}^*\phi{a\brack b}\left/\right.
\phi{a+e\brack b+f}=\exp(-2\pi i\ {}^teb).$$
Taking $N$th powers gives the result, 
when $d$ is as defined earlier in this section.
\end{proof}
\end{lemma}

On the other hand, construct for every $P$ a Weil function $f_P$
by our determinantal formula above. 
Then construct the $N$th, or $2N$th, powers
$(\tf_P^N)_{P\in A[N]}$, or $(\tf_P^{2N})_{P\in A[N]}$, of
a normal Weil set according to \ref{powers}.

\begin{theorem}\label{thomae} 
\part[i] There is a normal Weil set $(\tf_P)_{P\in A[N]}$ such that
$$\tf_P=\left(\theta{-a\brack -b}(z,\tau)\left/\right.\theta{0\brack 0}(z,\tau)\right)^{N}.$$

\part[ii] Choose $w=(w_1,\ldots,w_{3g-3})\in C^{3g-3}$
with $\sum w_i$ linearly equivalent to $K_C+\delta$. Then
$\tf_P(w)=(\theta{-a\brack -b}(0,\tau)/\theta{0\brack 0}(0,\tau))^N$.
\noproof
\end{theorem}

\begin{remark}\label{thetanull}
Let us recapitulate the situation.

\part[i] Starting with a curve $C$ of genus $g\ge 2$,
embed $X=\Jac^{g-1}_C$ into $\P^{3^g-1}$ using the polynomials
from \ref{polynomials}.

\part[ii] Given an integer $N$ prime to $\ch k$, compute the $N$-torsion
points on $A=\Jac^0_C$. For example, use the chord and tangent process described
in Section \ref{group law} or Section{computing torsion}.

\part[iii] For every $N$-torsion divisor class $D_P$ on $C$, 
compute a basis $\{\sigma_i^P\}$ of $H^0(C,\sO_C(2K_C+D_P))$.

\part[iv] Define the Weil function 
$f_P=\det(\sigma_i^P(z_j))^N/\det(\sigma_i^0(z_j))^N$.
Via the abelian sum $\alpha:C^{3g-3}\to X$, this is a rational function
on both $X$ and $C^{3g-3}$.

\part[v] Compute the Weil pairing $e_N$ from \ref{weil pairing}.

\part[vi] Choose an asymmetric bilinear pairing $d_N$ on $A[N]$
whose skew-symmetrization is $e_N$.

\part[vii] Compute the normalized Weil
functions $(\tf_P)_{P\in A[N]}$; their $2N$th powers are given 
by \ref{powers} (and their $N$th powers if $N$ is odd), 
and then finding $\tf_P$ is a matter of a finite search.

\part[viii] The functions $(\tf_P)_{P\in A[N]}$ define
a morphism $X\to P^{N^{2g}-1}$; the image of the subscheme $Fix_X$
is the moduli point of (the Jacobian of) $C$ with level structure
(level $N$ if $N$ is odd, level $(N,2N)$ if $N$ is even).

\part[ix] Suppose that $k=\C$ and that $N$ is even. A symplectic basis of
$H_1(C,\Z)$ determines a period matrix $\tau$ and a
level $(N,2N)$ structure, and then a level $2$ structure,
which then determines
an isomorphism $X\to A$, by \ref{igusa}. Our analogue of Thomae's
formula is
$$\tf_P(w_1,...,w_{3g-3})^{2N}=
(\theta{-a\brack -b}(z,\tau)/\theta{0\brack 0}(z,\tau))^{2N^2}$$
when the points $(w_1,...,w_{3g-3})$ of $C^{3g-3}$ and $z$ in the universal
cover $\C^g$ of $A$
have the same image in $A$, and $P=\tau a+b$.

Recall also that $\theta{-a\brack -b}(0,\tau)$, $\theta{-a\brack -b}(0,\tau)^2$ 
and $\theta{-a\brack -b}(0,\tau)^4$
are modular forms on $\Gamma_g(4,8)$, $\Gamma_g(2,4)$ and $\Gamma_g(2)$,
of weights $1/2$, $1$ and $2$.
\end{remark}

The disadvantage of this result is that it is not clear how to write down
a symplectic isomorphism $\phi$. However, the natural way of avoiding
this is to take invariants. If $N$ is even then the natural
finite group acting is $G=Sp_{2g}(\Z)/\Gamma(N,2N)$, which is an
extension of $Sp_{2g}(\Z/N)$ by $(\Z/2)^{2g}$, and the stack quotient
$[\sA_{g,(N,2N)}/G]$ is just $\sA_g$; we can pass to geometric quotients
to make things concrete.
\end{section}
\bigskip
\begin{section}{Spinorial square roots of theta-nulls on $\sM_g$}\label{spinor}
\medskip
Suppose for the moment that the base field is $\C$.
If ${a\brack b}$ is a half-integer characteristic, then the theta-null
$\theta{a\brack b}=\theta{a\brack b}(0,\tau)$ is a modular form of 
weight $\frac{1}{2}$ and level $(4,8)$.
These divisors are conjugate under the Galois group
$Sp_{2g}(\Z)/\Gamma_g(4,8)$, 
so they define a reduced and irreducible divisor $\theta_0$ on $\sA_g$
(because, for example, the second $\ell$-adic Betti number
of $\sA_g$ is $1$ and there are no modular forms of weight less than 
$\frac{1}{2}$). 
When $\sA_g$ is identified with the stack of principal symmetric
abelian torsors $(X,\Theta)$, then $\theta_0$ is the locus
where there is a point of $Fix_X$ that lies on $\Theta$
with positive even multiplicity. This locus is a divisor
on $\sA_g$ over $\Sp\Z[1/2]$.  

This section considers 
the local geometry of $\theta_0$. Unsurprisingly, this
reflects the geometry of the theta divisor $\Theta_\tau$
at the origin $0_\tau$ of the ppav $A_\tau$
corresponding to a point $\tau$ of $\theta_0$; in fact,
the heat equation shows, at once, that if the theta divisor
$\Theta_\tau$ on $A_\tau$
is of multiplicity $2n$ at the origin $0_\tau$ of $A_\tau$, 
then $\theta{0\brack 0}$ is of multiplicity exactly $n$ at $\tau$.
The determinantal formulae for theta functions on Jacobians
give something more detailed, on $\sM_g$, however.
To describe this needs a
well known elementary construction from representation theory,
which we now recall. Suppose now that the base field 
is of characteristic $p\ne 2$.

Suppose that $V$ is a $2n$-dimensional vector space with
a non-degenerate symmetric bilinear form $\phi$ and
associated quadratic form $q$. Fix a maximal isotropic
subspace $V_0$ of $V$. Then the set of maximal isotropic
subspaces $U$ of $V$ such that $\dim(U\cap V_0)$ is \emph{even}
is a spinor variety $S=Spin_n/P$; in particular, it is irreducible.
Forgetting the quadratic form gives an embedding
$i:S\inj G=Grass(n,2n)$ such that the relationship between
the tautological line bundles on them is
$i^*\sO_G(1)\cong\sO_S(2)$. We can regard $S$ as the variety of pure spinors
in $\P(W)$, where $W$ is a $2^{n-1}$-dimensional spin representation of $Spin_n$; then
$H^0(S,\sO_S(1))$ is the dual representation $W^\vee$. Moreover, the subspace $V_0$
defines hyperplane sections $H'$ and $H$ of $S$ and $G$,
which are closures of Schubert cells. Geometrically,
$H$ is the locus of $n$-dimensional subspaces of $V$ 
that meet $V_0$ non-trivially and $H'$ 
is the locus of isotropic $n$-dimensional subspaces of $V$
that meet $V_0$ non-trivially, so in a vector space of dimension at least two. 
That is, $H$
is defined by the vanishing of one Pl{\"u}cker co-ordinate $s$, say,
an $n\times n$ determinant that, given an $n$-dimensional
subspace $W$ of $V$, defines the degeneracy
locus of the composite map $W\to V/V_0$,
while $H'$ is defined by the vanishing of an element $v\in W^\vee$,
which is a pure spinor, such that 
$v=\sqrt{i^*s}$. In other words, $H$ has contact with $S$ along $H'$,
or $H.S=2H'$, where this is an equality of divisors, not just
a linear equivalence.

The first part of the next result, the existence
of a square root of a theta null on $M_g$, 
is due to Tsuyumine (\cite{Ts}, Theorem $1$).
That this square root is spinorial is the second part.

\begin{theorem}\label{square_root} 
\part[i] The restriction of $\theta_0$
to $\sM_g$ is a divisor $E$ of multiplicity $2$.
That is, $\theta_0$ has contact along $\sM_g$.

\part[ii] The square root $E$ of $D$ appearing in \DHrefpart{i}
is locally modelled by the embedding $i:S\inj G$
described above, with $n=3g-3$. That is, locally on $\sM_g$ there
is a morphism $\pi:\sM_g\to S$ with $E=\pi^*H'$ and
$D=\pi^*i^*H$. 

\begin{proof} We start by giving a short proof
of a crude version of \DHrefpart{i}, namely, that
$E$ has multiplicity at least $2$ everywhere.

By the Riemann--Kempf singularity theorem, if $C$ is a curve
whose moduli point lies on $\theta[0,0]$, then
$C$ has an effective even theta-characteristic $D$.
That is, $2D\sim K_C$, $h^0(C,\sO_C(D))$ is even
and at least $2$. Pick a point $x=\sum P_j\in C^{3g-3}$ that is
not on $\Delta$ and is linearly equivalent to $K_C+D$;
such an $x$ exists because the linear system
$\vert K_C+D\vert$ has no base points.
Choose a basis $\{\sigma_i\}$ of $H^0(C,\sO(2K_C))$.
Since $h^0(C,\sO(2K_C-x))=h^0(C,\sO(D))$, by Serre duality,
it follows that the corank of the matrix $(\sigma_i(P_j))$
is $h^0(C,\sO(D))$, which is even. Since $H^0(C,\sO(2K_C))$
is naturally isomorphic to the cotangent space of $\sM_g$
at the point $[C]$, what we have is a morphism $\phi$
from the germ $(\sM_g,[C])$ of $\sM_g$ at $[C]$ to the space
of $(3g-3)\times(3g-3)$ matrices such that $\phi([C])$
lies in the locus of matrices
whose corank is at least $h^0(C,\sO(D))$. That is,
$\phi$ maps $(\sM_g,[C])$ 
to the singular locus of the hypersurface $\det =0$.

To prove the full version we elaborate the ideas of \cite{MuTh}
so as to cover all line bundles of degree $g-1$, not only theta characteristics,
and to let the curve vary.

Start with a curve $C$ of genus $g$ over an algebraically closed field
$k$. Allow $\ch k$ to be arbitrary. Fix an integer $m\ge 2$.

Suppose that $E,F$ are line bundles
on $C$ of degree $g-1$, with a perfect bilinear pairing
$\Phi:E\times F\to \omega_C$.
Of course, $\Phi$ is an isomorphism on the tensor product
and $F\cong \omega_C\otimes E^\vee$.
Suppose that $a,b$ are effective divisors on $C$ of degree $(2m-1)(g-1)$;
then $\Phi$ extends to a perfect bilinear pairing
$\Phi:E(a)\otimes F(b)\to\omega(a+b)$.

Define $\phi: H^0(C,E(a)/E(-b))\times H^0(C,F(b)/F(-a))\to k$ by 
$$\phi({\bar e},{\bar f})=\sum_{P_i\in a+b} \Res_{P_i}\Phi(e_i,f_i),$$
where $e_i, f_i$ are liftings of ${\bar e}, {\bar f}$ to local sections 
of $E(a), F(b)$ and the sum is taken with multiplicities.

\begin{lemma}\label{symmetric} $\phi$ is a well defined perfect bilinear pairing.

\begin{proof} I omit the proof that $\phi$ is well defined and
bilinear. 

Suppose that ${\bar e}\in H^0(C,E(a)/E(-b))$
and is non-zero at the point $P_i$ in the support of $a$.
Then pick ${\bar f}$ to lie in the subspace $H^0(C,F/F(-a)$
of $H^0(C,F(b)/F(-a))$, non-zero at $P_i$, zero elsewhere.
Then $\Res_{P_i}\Phi(e_i,f_i)\ne 0$ while the other residues vanish.
So $\phi({\bar e},{\bar f})\ne 0$. So $\phi$ is non-degenerate
on the left; non-degeneracy on the right is proved similarly.
\end{proof}
\end{lemma}

Now suppose that $a,b$ are chosen so that $E(a)$ is isomorphic
to $F(b)$, which is in turn isomorphic to $E^\vee\otimes\omega(b)$. 
Choose such isomorphisms and use them to identify these sheaves.
Twisting by $\sO_C(-a-b)$
then induces an identification $H^0(C,E(a)/E(-b))= H^0(C,F(b)/F(-a))=V,$
say.

\begin{lemma} With this identification, $\phi$ is symmetric.
\noproof
\end{lemma}

Abbreviate $(2m-1)(g-1)$ to $n$. Then $\dim V=2n$ and
the cohomology of the exact sequences
$$0\to E(-b)\to E(a)\to E(a)/E(-b)\to 0,$$
$$0\to E/E(-b)\to E(a)/E(-b)\to E(a)/E\to 0$$
exhibits $H^0(C,E(a))$ and $H^0(C,E/E(-b))$
as $n$-dimensional subspaces of the quadratic space $(V,\phi)$.

\begin{lemma}\label{10.3} $H^0(C,E(a))$ is totally isotropic.

\begin{proof} Suppose that ${\bar e},{\bar f}\in H^0(C,E(a))$.
Then in the definition of $\phi({\bar e},{\bar f})$
we can take $e_i={\bar e}$ and $f_i={\bar f}$, so that
$\phi({\bar e},{\bar f}) =\sum \Res_{P_i}\Phi({\bar e},{\bar f})$.
This vanishes, by the residue theorem, since $\Phi({\bar e},{\bar f})$ 
is a global rational $1$-form on $C$.
\end{proof}
\end{lemma}

So $\phi$ induces an exact sequence
$$0\to H^0(C,E(a))\to V\to H^0(C,E(a))^\vee\to 0.\ \ (*)$$
The two exact sequences preceding \ref{10.3} fit into
an exact commutative diagram of sheaves on $C$, 
as follows:
$$\xymatrix{
&&{0}\ar[d] & {0}\ar[d]\\
{0}\ar[r]&{E(-b)}\ar[d]_{=}\ar[r]&{E}\ar[d]\ar[r]&{E/E(-b)}\ar[d]\ar[r]&{0}\\
{0}\ar[r]&{E(-b)}\ar[r]&{E(a)}\ar[d]\ar[r]&{E(a)/E(-b)}\ar[d]\ar[r]&{0}\\
&&{E(a)/E}\ar[r]^{=}\ar[d]&{E(a)/E}\ar[d]\\
&&{0} &{0}
}$$
The cohomology of this and the 
duality isomorphism
$$H^1(C,E(-b))\cong H^0(C,E^\vee\otimes\omega(b))^\vee
= H^0(C,F(b))^\vee= H^0(C,E(a))^\vee$$
give a commutative exact diagram
$$\xymatrix{
&{0}\ar[d] & {0}\ar[d]\\
{0}\ar[r]&{H^0(C,E)}\ar[d]\ar[r]&{H^0(C,E/E(-b))}\ar[r]\ar[d]\ar[dr]^{\xi}& 
{H^0(C,E(a))^\vee}\ar[d]^{=}\\
{0}\ar[r]&{H^0(C,E(a))}\ar[d]_{\mu}\ar[r]&{V}\ar[r]\ar[d]^{\nu}&
{H^0(C,E(a))^\vee}\ar[r]&{0}\\
&{H^0(C,E(a)/E))}\ar[r]^{=}&{H^0(C,E(a)/E))}\ar[d]\\
&&{0}
}$$
The middle row
of this diagram is identical to the exact sequence $(*)$ above.

Denote the subspaces $H^0(C,E(a))$ and $H^0(C,E/E(-b))$ of $V$
by $V_0$ and $V_1$ respectively, and put $H^0(C,E(a)/E)=V_2$,
so that $V_2$ is identified with $V/V_1$.

\begin{lemma} In this diagram, 
$H^0(C,E)=H^0(C,E(a))\cap H^0(C,E/E(-b))$ inside $V$.
That is, $\ker\xi= H^0(C,E)$.

\begin{proof} Certainly, $H^0(C,E)\subset H^0(C,E(a))\cap H^0(C,E/E(-b))$.
Conversely, suppose $x\in H^0(C,E(a))\cap H^0(C,E/E(-b))$. Then
$\nu(x)=0$, so $\mu(x)=0$, so $x\in H^0(C,E(a))\cap H^0(C,E/E(-b))$.
\end{proof}
\end{lemma}

Next, suppose that $E\in J=\Jac^{g-1}_C$ and $a\in C^{(n)}$ are allowed to
vary, but subject to the constraint that $E(a)\cong\omega_C^{\otimes m}$.

Take $F,b$ as before; this determines $F$ up to isomorphism
and $b$ up to linear equivalence. In particular, $b\sim (m-1)K_C+E$.
These constraints define subvarieties
$Z$ of $J\times C^{(n)}\times C^{(n)}$
and $Z_1,Z_2$ of $J\times C^{(n)}$; that is,
$Z=\{(E,a,b)\vert E(a)\cong\omega_C^{\otimes m},\ b\sim (m-1)K_C+E\}$,
$Z_1=\{(E,b)\vert E(-b)\cong\omega_C^{1-m}\}$ and
$Z_2=\{(E,a)\vert E(a)\cong\omega_C^{\otimes m}\}$.
These fit into a commutative square
$$\xymatrix{
{Z}\ar[r]^{\delta=pr_{12}}\ar[d]_{\beta =pr_{13}}\ar[dr]^{\alpha}& 
{Z_2}\ar[d]^{\epsilon =pr_1}\\
{Z_1}\ar[r]_{\gamma=pr_1}& {J}
}$$
where $\beta,\gamma,\delta,\epsilon$ are,
by Abel's theorem, $\P^{n-g}$-bundles
and $\alpha:Z\to J:(E,a,b)\mapsto E$ is 
a $\P^{n-g}\times \P^{n-g}$-bundle.
In fact, this square is, obviously, Cartesian.
Moreover, there is a universal line bundle $E$ over $Z\times C$ and universal
cycles $a,b\subset Z\times C$. 

Let $pr_Z:Z\times C\to Z$ be the projection 
and put $\sV=pr_{Z*}E(a)/E(-b)$,
$\sV_0=pr_{Z*}E(a)$, $\sV_1=pr_{Z*}E/E(-b)$
and $\sV_2=pr_{Z*}E(a)/E$; 
these are locally free sheaves
of ranks $2n$, $n$, $n$ and $n$ respectively,
and $\sV_0\cong H^0(C,\omega_C^{\otimes m})\otimes\sO_Z$. 
Exactly as before, $\sV$ has a non-degenerate
symmetric bilinear pairing $\phi:\sV\times\sV\to\sO_Z$ and $\sV_0$
is identified with a maximal isotropic 
sub-bundle of $\sV$. In particular, $\phi$ induces an exact sequence
$$0\to\sV_0\to\sV\to\sV_0^\vee\to 0$$
and there is another exact sequence
$$0\to\sV_1\to\sV\to\sV_2\to 0.$$

\begin{proposition}\label{vary}\label{pullpush}
\part[i]The triple $(\sV,\phi,\sV_0)$ pulls back via $\alpha$
from a triple
$(\sU,\psi,\sU_0)$ on $J$, where $\psi$ is a non-degenerate
symmetric bilinear form on $\sU$ and $\sU_0$
is a maximal isotropic sub-bundle.
 
\part[ii] $0\to\sV_1\to\sV\to\sV_2\to 0$ pulls back from an exact 
sequence $0\to\sU_1\to\sU\to\sU_2\to 0$ of sheaves on $J$.

\begin{proof} Define $\sU_0=H^0(C,\omega_C^{\otimes m})\otimes\sO_J$.
Then $\sV$ is a class in $\Ext^1_Z(\alpha^*\sU_0,\alpha^*\sU_0^\vee)$.
Since $\alpha$ is a $\P^{n-g}\times\P^{n-g}$-bundle,
it follows that $R^1\alpha_*\sO_Z=0$, and so
$$\Ext^1_Z(\alpha^*\sU_0,\alpha^*\sU_0^\vee)=\Ext^1_J(\sU_0,\sU_0^\vee).$$
So $\sV=\alpha^*\sU$ for some $\sU\in \Ext^1_J(\sU_0,\sU_0^\vee)$.

The existence of the quadratic form $\psi$ follows similarly:
$\Hom_Z(\Symm^2\sV,\sO_Z)=\Hom_J(\Symm^2\sU,\sO_J)$.

The proof of \DHrefpart{ii} is also similar in outline.
Since the morphisms $\alpha,\ldots,\epsilon$
in the Cartesian square above
are proper and flat, 
there is, by the flat base change theorem for coherent cohomology,
an equivalence
$\gamma^*R^q\epsilon_*=R^q\beta_*\delta^*$ of functors
on the category of coherent sheaves $\sF$ on $Z$
that are flat over $J$.

By construction, there are locally free sheaves $\sG_1$ and $\sG_2$ on 
$Z_1,Z_2$ respectively such that $\sV_1=\beta^*\sG_1$
and $\sV_2=\delta^*\sG_2$. So $0\to\sV_1\to\sV\to\sV_2\to 0$
equals $0\to\beta^*\sG_1\to\alpha^*\sU\to\delta^*\sG_2\to 0$.
Define $\sW_1=\beta_*\delta^*\sG_2$; then the vanishing of
$R^1\beta_*\sO$ gives an exact sequence
$0\to\sG_1\to\gamma^*\sU\to\sW_1\to 0$ on $Z_1$, and it follows
at once that $\sW_1$ is locally free and $\beta^*\sW_1=\delta^*\sG_2$.
The same argument applied to the dual sequence gives a locally
free sheaf $\sW_2$ on $Z_2$ with $\delta^*\sW_2=\beta^*\sG_1$.

The equivalence of functors above shows that
$$\gamma^*R^q\epsilon_*\sG_2=R^q\beta_*\delta^*\sG_2
=R^q\beta_*\beta^*\sW_1=\sW_1\otimes R^q\beta_*\sO_Z,$$
which vanishes when $q\ge 1$ and equals $\sW_1$ when $q=0$.
By faithful flatness, this gives $R^q\epsilon_*\sG_2=0$
for $q\ge 1$.

Similarly, $R^q\gamma_*\sW_1=0$ for $q\ge 1$
and $\epsilon^*\gamma_*\sW_1=\sG_2$.

So take $\sU_1=\gamma_*\sW_1$ and $\sU_2=\epsilon_*\sW_2$.
That these sheaves have the properties we seek is now
immediate by the argument used before involving $Ext^1$.
\end{proof}
\end{proposition}

By abuse of notation, we also let $\xi$ denote the
composite homomorphism $\sU_1\to\sU\to\sU_0^\vee$
of locally free sheaves on $J$.

In crude terms, we have shown that
the quadratic vector space $(V,\phi)$,
its maximal isotropic subspace $V_0$, its subspace $V_1$ 
and its quotient space $V_2=V/V_1$
depend upon $E$ but are independent of $a$ and $b$ when $E,a$ and $b$
are permitted to vary subject to our constraint.
In particular, the composite homomorphism
$\xi:V_1\to V/V_0=V_0^\vee$
depends only on $E$.

The conclusion to be drawn is that, if $\Theta_J\subset J$
is the theta divisor, then 
$\Theta_J$ is exactly the locus where $\xi$ is degenerate. Since the target
of $\xi$ is $H^0(C,\sO_C(mK_C))^\vee$, this coincides with
the previous determinantal description of $\Theta$
after pulling back to $C^{n}$ under the abelian sum.

\begin{lemma}\label{theta char}\label{isotropic} 
Suppose that $E$ is a theta-characteristic. Then
$V_1$ is also an isotropic subspace of $V$.
\begin{proof}
This is \cite{MuTh}, p. 184, first paragraph.
\end{proof}
\end{lemma} 

Now relativize this construction; that is, let the curve $C$,
the line bundle $E$ of degree $g-1$
and the effective cycles $a,b$ of degree $n$
all vary, subject to the same constraints
that $E(a)\cong\omega_C^{\otimes m}$ and
$\sO_C(b)\cong E \otimes\omega_C^{\otimes(m-1)}$.

So suppose that $f:C\to S$ is a family
of genus $g$ curves; this leads to morphisms $\beta:J=\Jac^{g-1}_{C/S}\to S$
and $C^{(n)}\to S$. There is a closed subscheme $Z$ of
$J\times_S C^{(n)}\times_S C^{(n)}$ whose points consist
of triples $(E,a,b)$ as above
The projection $\pi:Z\to S$ factors as
$Z\stackrel{\alpha}{\to}J\stackrel{\beta}{\to} S$, where
$\alpha$ is a $\P^{n-g}\times\P^{n-g}$-bundle.

The quadratic vector space $(V,\phi)$ is replaced by a
quadratic vector bundle $(\sV,\phi)$ over $Z$, and the maximal
isotropic subspace $V_0$ by a sub-bundle $\sV_0$ which is isomorphic
to $\pi^*\sF_m$, where $\sF_m=f_*\omega_{C/S}^{\otimes m}$. So
$\sV$ appears as a class in
$\Ext^1_Z(\pi^*\sF_m^\vee,\pi^*\sF_m)$.

\begin{proposition} 
\part[i] The
triple $(\sV,\phi,\sV_0)$ is the pullback under $\alpha:Z\to J$
of some triple $(\sU,\psi,\sU_0=\beta^*\sF_m)$ on $J$
where $\sU$ defines a class
$[\sU]$ in 
$\Ext^1_J(\beta^*\sF_m^\vee,\beta^*\sF_m)$.

\part[ii]
The exact sequence $0\to\sV_1\to\sV\to\sV_2\to 0$
is the pullback under $\alpha$ of an exact sequence
$0\to\sU_1\to\sU\to\sU_2\to 0$ on $J$.

\part[iii] The degeneracy divisor of the composite homomorphism
$\xi:\sU_1\to\sU_0^\vee$ is the theta divisor $\Theta$ on $J$.
\begin{proof} Exactly as \ref{vary}
\end{proof}
\end{proposition}

This can be rephrased as a theorem about $\sM_g$
over $\Z$.

Let $f:C\to\sM_g$ be the universal curve,
$\sF_m=f_*\omega_{C/\sM_g}^{\otimes m}$, 
$\beta:J=\Jac_C^{g-1}\to\sM_g$
its degree $g-1$ Jacobian and $\Theta=\Theta_J\subset J$
the theta divisor. We have taken $m\ge 2$ and put
$n=(2m-1)(g-1)$.

\begin{theorem}\label{modular bundles} 
\part[i] There is a rank $2n$ vector bundle $\sU$ on $J$ 
with a non-degenerate symmetric bilinear form $\phi$ and a maximal
isotropic sub-bundle $\sU_0$ isomorphic to $\beta^*\sF_m$.

\part[ii] There is a rank $n$ sub-bundle $\sU_1$ of $\sU$ such 
that the degeneracy locus of the composite homomorphism $\xi:\sU_1\to\sU_0^\vee$
is $\Theta$.
\noproof
\end{theorem}

Now exclude characteristic $2$. That is, work over $\Z[\frac{1}{2}]$.
We shall work at level $(2,4)$. 

Take the universal principal symmetric torsor
$(X,\Theta)\to \sA_{g,(2,4)}$ over the stack $\sA_{g,(2,4)}$ of ppav's with 
level $(2,4)$ structure and automorphism group scheme $U=\Aut_X^0\to\sA_{g,(2,4)}$.
As before, the induced level $2$ structure determines an even theta characteristic
$\Sigma$, 
so an isomorphism $U\to X$ taking the zero-section $0_U$ of $U\to\sA_{g,(2,4)}$ to
the section $\Sigma$ of $X\to \sA_{g,(2,4)}$.

According to what was recalled in section \ref{prelim},
the level $(2,4)$ structure determines the
normalized Weil functions $(\tf_P)_{P\in U[2]}$.
Over $\C$, these can be expressed as ratios of
squares of theta functions $\theta[a,b](z,\tau)$ with half-integer characteristics.

The functions $\tf_P$ define the morphism
$U\to\P^{2^{2g}-1}$ recalled in section \ref{prelim}. In particular,
they give the canonical reference hyperplanes in $\P^{2^{2g}-1}$,
and, by construction, each of these hyperplanes cuts $U$ with
multiplicity $2$.
Evaluating these functions at $\Sigma$ is the morphism
$\pi:\sA_{g,(2,4)}\to\P^{2^{2g}-1}$; again, each of the canonical hyperplanes
cuts $\sA_{g,(2,4)}$ along the divisor $2\theta_P$ of multiplicity $2$, 
where $\theta_P$ is empty if $P$ is an odd $2$-torsion point
and $\sum_{P\ odd}\theta_P$ is the pullback to $\sA_{g,(2,4)}$
of the divisor $\theta_0$ on $\sA_g$.
Of course, over $\C$, this morphism $\pi$ is defined by the squares 
$\theta[a,b](0,\tau)^2$ of the theta-constants and there is an equality
$\theta_P=(\theta[a,b](0,\tau))_0$ of divisors on $\sA_{g,(2,4)}$
where $P=\tau a+b$. The thetanull $\theta[a,b](0,\tau)$ is  
a modular form (of weight $\frac{1}{2}$) on the group $\Gamma(4,8)$,
not on $\Gamma(2,4)$ (on $\Gamma(2,4)$ it has a quadratic character), 
but its square $\theta[a,b](0,\tau)^2$
is a modular form on $\Gamma(2,4)$, without any character.

Now pull back $(X,\Theta)$ via the Jacobian morphism $j:\sM_{g,(2,4)}\to \sA_{g,(2,4)}$. 
We get $\beta:(J,\Theta_J)\to \sM_{g,(2,4)}$, the degree $g-1$ Jacobian of the universal
curve $f:C\to \sM_{g,(2,4)}$. The section $\Sigma$ of
$X\to\sA_{g,(2,4)}$ pulls back to
a section $\Gamma$ of $J\to\sM_{g,(2,4)}$. Identify $\Gamma$ with $\sM_{g,(2,4)}$.
Consider the morphism $\Gamma\to\sA_{g,(2,4)}$ (which is not quite an embedding
because of the failure of local Torelli for hyperelliptic curves).
Then $\Gamma$ meets $J_t$ in a point corresponding to 
an even theta-characteristic $E_t$ on $C_t$, and
$\Gamma$ meets $\Theta_{J_t}$ if and only if $E_t$ is effective.
That is, $\Gamma$ meets $\theta_0$ along the locus of points
$t\in \Gamma=\sM_{g,(2,4)}$ where the curve $C_t$
has an effective even theta-characteristic,

Write $\sW,\sW_i$ for the restrictions of the bundles $\sU,\sU_i$ of
\ref{modular bundles} to $\sM_{g,(2,4)}$.

\begin{theorem} $\sW$ is a rank $2n$ vector bundle on $\sM_{g,(2,4)}$,
with a non-degenerate symmetric bilinear
form $\psi$. There are two maximal isotropic sub-bundles $\sW_0$
and $\sW_1$, with $\sW_0\cong\sF_m$, such that
the determinant of the induced homomorphism $\xi:\sW_1\to\sW_0^\vee$,
which is the zero locus of a theta-null on $\sM_{g,(2,4)}$, 
is the square of a pure spinor. 

\begin{proof} This follows at once from what we have done.
In particular, that $\sW_1$ is also isotropic is \ref{theta char}
and the existence of the spinorial square root is the
piece of representation theory at the start.
\end{proof}
\end{theorem}

Theorem 10.1 is a restatement of this.
\end{proof}
\end{theorem}

\begin{remark}
\part[i] The non-degenerate bilinear form $\phi$ on
the vector space $V$ is exactly what is given by the duality theorem
applied to the family of complexes $\sC^\bullet$ supported in degrees $-1,0$
with $\sC^{-1}=\omega_C^{\otimes 1-m}$, $\sC^0=\omega_C^{\otimes m}$
and whose differential is parametrized by the variety 
$Z\subset \Jac^{g-1}\times C^{(n)}\times C^{(n)}$ above.
Note that $\sC^\bullet\cong\sC^{\bullet\vee}\otimes\omega_C^\bullet$,
where $\omega_C^\bullet\cong\omega_C[1]$
is the dualizing complex.
Can this picture be embedded in a bigger one that
leads to some structure involving multiplication, or convolution,
of the bundles $\sV_m$ on $J$, and does this in turn lead to
any other structure on $\sM_g$?

\part[ii] When $g=3$, the Torelli morphism is simply ramified
along the hyperelliptic locus $\sH_3$, which is a divisor in $\sM_3$,
and $\theta[0,0]$ defines $\sH_3$ with multiplicity $2$.
When $g=2$, the theta-constants have no zeroes on $\sM_2$; they
vanish in $\sA_2$ exactly along the locus of decomposable abelian surfaces.
\end{remark}
\end{section}
\bigskip
\begin{section}{Formulae for Prym varieties}\label{prym}
\medskip
We start by considering various stacks of non-abelian bundles on $C$.

Denote by $GL_{2,C,K_C}$ the stack of rank $2$ vector bundles $\sH$ on $C$
with determinant $\omega_C$ and by $GL_{2,C,K_C,m}$,
for any integer $m\ge 1$, the open substack
defined by the conditions that $H^0(C,\sH(-mK_C))=0$ and
$H^1(C,\sH(mK_C))=0$. For any such $\sH$, put $\sH'=\sH^\vee\otimes\omega_C$,
Then, for all $a,b\in\vert mK_C\vert$,
there is an isomorphism $\sH(a)\to\sH'(b)$, unique up to scalars,
and so, via the determinant, a non-degenerate skew-symmetric
pairing $\Psi:\sH(a)\otimes\sH(a)\to\omega_C(a+b)$.
Define $V=H^0(C,\sH(a)/\sH(-b))$; its dimension is $8m(g-1)$.

\begin{lemma} There is a non-degenerate skew-symmetric bilinear pairing
$\psi:V\times V\to k$ defined by
$$\psi({\bar{e}},{\bar{f}})=\sum_{P_i\in a+b}\Res_{P_i}\psi(e_i,f_i).$$
\begin{proof}
Same as in the symmetric case, \ref{symmetric}.
\end{proof}
\end{lemma}

Suppose that $\sH$ lies in $GL_{2,C,K_C,m}$. Then, just as before,
if we define $V_0=H^0(C,\sH(a))$, $V_1=H^0(C,\sH/\sH(-b))$
and $V_2=H^0(C, \sH(a)/\sH)$, then
$V_0$ is a maximal totally isotropic subspace and
there are short exact sequences
$$0\to V_0\to V\to V_0^\vee\to 0,$$
$$0\to V_1\to V\to V_2\to 0.$$

Let $\vert mK_C\vert_1$ be the copy of $\vert mK_C\vert $ in which $a$ moves
and $\vert mK_C\vert_2$ that in which $b$ moves. 
Put $X=GL_{2,C,K_C,m}\times\vert mK_C\vert_1\times\vert mK_C\vert_2$,
$Y=GL_{2,C,K_C,m}\times\vert mK_C\vert_1$, 
$Z=GL_{2,C,K_C,m}\times\vert mK_C\vert_2$ and consider the
$2$-Cartesian square
$$\xymatrix{
{X}\ar[r]^q\ar[d]_p & {Z}\ar[d]\\
{Y}\ar[r] & {GL_{2,C,K_C,m}}.
}$$
By construction, the vector spaces $V,V_0,V_1,V_2$
and the pairing $\psi$ globalize to vector bundles
$\sV$ etc. and a pairing $\psi$ on $X$.
Moreover, there are vector bundles $\sU_0',\sU_2'$ on $Y$
such that $\sV_0=p^*\sU_0'$, $\sV_2=p^*\sU_2'$ and a vector
bundle $\sW_1'$ on $Z$ with $\sV_1=q^*\sW_1'$.

\begin{proposition} 
\part[i] The vector bundles $\sV,...,\sV_2$, the maps
between them and
the pairing $\psi$ all pull back from vector bundles
$\sU,...,\sU_2$, maps and a pairing $\psi$ on $GL_{2,C,K_C,m}$.
\part[ii] The theta divisor on $GL_{2,C,K_C,m}$ is the
degeneracy locus of the composite homomorphism 
$\sU_1\to\sU_0^\vee$. 
\begin{proof}
Same as \ref{pullpush}
\end{proof}
\end{proposition}

\begin{remark} In the abelian case, the bundle $\sV_0$ was,
when $m\ge 2$, 
the trivial bundle $H^0(C,\omega_C^{\otimes m})$.
It is not clear whether $\sV_0$ is trivial here, however.
\end{remark}

Now consider the stack $SO_{2,C,K_C}$ of rank $2$ vector bundles
$\sH$ on $C$ with a non-degenerate \emph{symmetric}
bilinear form $\Phi:\sH\otimes\sH\to\omega_C$ such that
$\det\sH$ is isomorphic to $\omega_C$. (of course, the existence of
$\Phi$ implies that $(\det\sH)^{\otimes 2}\cong\omega_C^{\otimes 2}$.)
Define the open substack $SO_{2,C,K_C,m}$ in the same way as above.
Then for any $a,b\in\vert mK_C\vert$ and $\sH\in SO_{2,C,K_C}$,
the vector space $V$ as above has a non-degenerate symmetric 
bilinear form $\phi$. 

\begin{proposition} Suppose that
$\tau:SO_{2,C,K_C,m}\to GL_{2,C,K_C,m}$
is the forgetful morphism.
\part[i] Then $\tau^*\sU$
has a non-degenerate symmetric bilinear form with respect
to which the bundles $\tau^*\sU_0$ and $\tau^*\sU_1$
are maximal isotropic.
\part[ii] The theta divisor $\Theta_{GL_2}$ on $GL_{2,C,K_C,m}$
pulls back to twice the theta divisor $\Theta_{SO_2}$ on $SO_{2,C,K_C,m}$.
\part[iii] \label{SO_2} $\Theta_{SO_2}$ has a spinorial structure.
\begin{proof} As before, especially \ref{isotropic} for the fact that
$\tau^*\sU_1$ is isotropic.
\end{proof}
\end{proposition}

\begin{remark} As in the previous section,
this construction can be extended to one on the stack $\sM_g$.
In terms of pairs of reductive groups relative
to the stack $\sM_g$, this comes to the consideration of
$1\subset \GG_m$ in the
previous section and $SO_2\subset GL_2$ in this one. 
\end{remark}

Suppose that $\pi:\tC\to C$ is an {\'e}tale double cover of curves,
with $g(G)=g$ and $g(\tC)=\tg$, so that $\tg =2g-1$
and $C=\tC/\langle\iota\rangle$, where $\iota$ is an involution
of $\tC$ with no fixed points. Consider the norm
map $\Nm\tX:=:\Jac_{\tC}^{\tg-1}\to\Jac_C^{2g-2}$. Then 
(\cite{MuTh}, \cite{MuP}) the fibre
$\Nm^{-1}(K_C)$ is the locus in $\tX$ where
$\iota$ is equal to $[-1_{\tA}]$ (here $\tA=\Aut_{\tX}^0$).
It has two connected components, $P^+$ and $P^-$, and
each is a torsor under the Prym variety $P=\Aut_{P^+}^0$. 
As subvarieties of
$\tX:=\Jac_{\tC}^{\tg-1}$, $P^+$ is the set of degree $\tg-1$
line bundles $\sL$ on $\tC$ such that $h^0(\sL)$ is even,
while $P^-$ is the locus where $h^0(\sL)$ is odd. The theta
divisor $\tth$ on $\tX$ contains $P^-$, while
$\tth.P^+=2\Xi$ for a divisor $\Xi=\Xi_0$ on $P^+$
such that $(P^+,\Xi)$ is a principal symmetric abelian torsor.
Note that $P$ is identified with the subgroup of
$\tA:=\Jac^0_{\tC}=\Aut_{\tX}^0$ consisting of points $x$
such that $t_x$, the translation 
by $x$, preserves $P^+$.

Note that
the norm defines a morphism $\rho:P^+\to SO_{2,C,K_C,1}$.

\begin{theorem} The divisor $\Xi$ on $P^+$ has a 
spinorial structure.

\begin{proof} Pull back the spinorial structure of \ref{SO_2}
under $\rho$.
\end{proof}
\end{theorem}

Finally, we describe Prym theta functions in
terms of determinants.

Suppose that $M\ge 1$ is an integer.
For every point $x\in P[2M]$, $Mt_x^*\tth.P=2M\t_x^*\Xi$
and is linearly equivalent to $2M\Xi$. Regard $x$ as a point in
$\tA[2M]$ and take a basis $\{\tsigma^x_i\}$ of
$H^0(\tC,\sO(2K_{\tC+x}))$. 

\begin{proposition} The ratio
$$(\det((\tsigma^x_i(z_j))/\det((\tsigma^0_i(z_j)))^M$$
restricts to a Weil function on $P^+$ for the
point $x$ of order $2M$.

\begin{remark} More accurately, there is a Weil
function for $x$ on $P^+$ whose pull-back to
the inverse
image of $P^+$ in $\tC^{3(\tg -1)}$
is given by this expression. 
\end{remark}

\begin{proof} Exactly as for Weil functions
on $Jac_C^{3g-3}$.
\end{proof}
\end{proposition}

\end{section}
\bigskip
\begin{acknowledgements}{Support by GNSAGA at the University of Rome, La Sapienza, 
and the Institut Mittag-Leffler (Djursholm) is gratefully acknowledged.}
\end{acknowledgements}
\bibliography{abbrevs,alggeom,combinat,algebra,ekedahl}
\bibliographystyle{pretex}

\end{document}